\documentclass[11pt]{article}
\usepackage[T1]{fontenc}
\usepackage[utf8]{inputenc}
\usepackage[margin=1.15in]{geometry}
\usepackage{amsmath,amssymb,amsthm,bbm, mathtools}
\usepackage{enumitem}
\usepackage{xcolor}
\usepackage[hidelinks,bookmarksnumbered]{hyperref}

\numberwithin{equation}{section}
\theoremstyle{plain}
\newtheorem{theorem}{Theorem}[section]
\newtheorem{corollary}[theorem]{Corollary}
\newtheorem{proposition}[theorem]{Proposition}
\newtheorem{lemma}[theorem]{Lemma}
\newtheorem{fact}[theorem]{Fact}
\newcommand{\restatedname}{}
\newtheorem*{restated}{\restatedname}
\newenvironment{restate}[2]{\renewcommand{\restatedname}{#1~\ref{#2}}\begin{restated}}{\end{restated}}
\theoremstyle{definition}

\newtheorem*{definition*}{Definition}
\newtheorem{remark}[theorem]{Remark}
\setlist[enumerate]{label=(\arabic*)}
\title{The spectrum of free projective geometries over finite chain rings}
\author{Yishai Lavi}
\date{}

\begin{document}
\maketitle

\begin{abstract}
Let $\mathcal{O}$ be a Dedekind domain, $\mathfrak{p}$ a prime ideal of
residue order $q$, and $\mathcal{O}_r=\mathcal{O}/\mathfrak{p}^r$ which is a finite chain ring; the motivating example
is $\mathbb{Z}/p^r$. The free projective geometry $\mathbb{P}^{d-1}_{fr}(\mathcal{O}_r)$ is
the graph whose vertices are the nontrivial free submodules of $\mathcal{O}_r^{d}$,
two of them adjacent when one strictly contains the other. For every $d$, $q$ and $r$
we reduce its adjacency matrix to explicit symmetric matrices (blocks) of size
at most $d-1$. There is one block for each signature (a partition with at most $d/2$ parts, each at most
$r$) and we show that every block, at every depth, is a single universal matrix
$L(s,\sigma)$ evaluated at the number $s$ of parts of the signature $\lambda$ and
at $\sigma=2rs-|\lambda|$. The same blocks give the spectrum of the subgraph induced
on any set of levels, and for two levels this spectrum is completely explicit. All of
it depends on the ring only through $q$ and $r$; in particular the subgraphs
$\mathbb{P}^{d-1}_{S}$ over $\mathbb{Z}/p^r$ and over $\mathbb{F}_p[t]/(t^r)$ are
isospectral for all $\varnothing\ne S\subseteq\{1,\dots,d-1\}$; together with known
rigidity results of arXiv:2510.09324, this proves that they are isospectral and non-isomorphic whenever $S$ contains
at least two levels, as conjectured there.
\end{abstract}

\section{Introduction}\label{sec:1}

The spherical building of $PGL_d(\mathbb{F}_q)$ is the simplicial complex of
nontrivial subspaces of $\mathbb{F}_q^{\,d}$, whose cells are flags. It appears as the
link of a vertex in the Bruhat--Tits building of $PGL_d(F)$, for a
non-archimedean local field $F$ with residue field $\mathbb{F}_q$, and its
expansion - the fact that its nontrivial eigenvalues are small compared with the
trivial one - underlies Garland's method and its descendants in the theory of
high-dimensional expanders \cite{Gar73,Opp18}. When the vertex link is replaced by
a \emph{geodesic sphere} of radius $r$ in the building (see \cite{KP22}), one obtains instead the complex of
free submodules of $\mathcal{O}_r^{\,d}$, where $\mathcal{O}_r=\mathcal{O}/\mathfrak{p}^r$
is a finite quotient of a Dedekind domain, such as $\mathbb{Z}/p^r$ or
$\mathbb{F}_q[t]/(t^r)$ \cite{KP22}. These complexes, called \emph{free projective
geometries}, were studied in \cite{LP25}, which establishes two forms of rigidity for them and
computes the spectrum of the bipartite graph $\mathcal{X}_{1,n}$ induced by the free
modules of ranks $1$ and $n$, for $n\le d-1$. The computation shows that $\mathcal{X}_{1,n}$ is an
excellent expander, and that its spectrum depends on the ring only through its length
$r$ and its residue order $q$; combined with rigidity, this produces pairs of isospectral
non-isomorphic graphs over $\mathbb{Z}/p^r$ and $\mathbb{F}_p[t]/(t^r)$. It was
conjectured there that the same holds for every two-level graph $\mathcal{X}_{m,n}$
\cite[Conj. 6.2]{LP25}.

Our main results (Theorems~\ref{thm:4.4}, \ref{thm:6.15} and \ref{thm:7.6} below) prove a
stronger form of the conjecture. While the methods of \cite{LP25} were purely combinatorial,
our methods here involve the representations of $GL_d(\mathcal{O}_r)$.
In this paper we reduce the problem of finding the spectrum of the \emph{whole} geometry - all levels at once,
for every $d$, $q$ and $r$ - to finding the spectrum of explicit symmetric blocks of size
at most $d-1$, with the multiplicities in a closed form. This suffices to deduce isospectrality
(depending only on $r$, $q$ and $d$). For graphs of two levels, we can even give
explicitly the full spectrum. The eigenvalues of the blocks
themselves are not given in closed form in general.

For $r=1$, where
the geometry is the inclusion graph of the subspaces of
$\mathbb{F}_q^{\,d}$, our explicit description of the blocks can be seen as a simplification of the blocks described in \cite{Lew24}, whose eigenvalues
were analysed only asymptotically as $q\to\infty$.

\subsection{The setting}\label{sec:1.1}

Throughout, $\mathcal{O}$ is a Dedekind domain, $\mathfrak{p}$ a nonzero prime ideal
of residue order $q=|\mathcal{O}/\mathfrak{p}|$, $r\ge1$,
$\mathcal{O}_r=\mathcal{O}/\mathfrak{p}^r$, and $d\ge3$. The main example is
$\mathcal{O}_r=\mathbb{Z}/p^r$; others are $\mathbb{F}_q[t]/(t^r)$ and
$\mathcal{O}_K/\mathfrak{p}^r$ for a number field $K$. The ring $\mathcal{O}_r$ is
a finite chain ring, of size $q^r$. An
\textbf{$n$-space} is a free submodule of $V=\mathcal{O}_r^{\,d}$ of rank $n$ which is a direct
summand, and $\mathcal{X}_n$ denotes the set of $n$-spaces, for $1\le n\le d-1$. We call $\mathcal{X}_n$ a
\textbf{level}, and, as in \cite{LP25}, we also say that its vertices have \textbf{color} $n$; we use
both terms. The free projective geometry $\mathbb{P}^{d-1}_{fr}(\mathcal{O}_r)$ is the graph with
vertex set $\mathcal{X}=\bigsqcup_n\mathcal{X}_n$ in which two spaces are adjacent when one strictly
contains the other, and $A$ is its adjacency matrix. For a set $S$ of levels we write
$\mathcal{X}_S$ for the induced subgraph on $\bigsqcup_{n\in S}\mathcal{X}_n$, and
$\mathcal{X}_{m,n}=\mathcal{X}_{\{m,n\}}$. The group $G=GL_d(\mathcal{O}_r)$ acts on
$\mathcal{X}$ preserving inclusion. We write $[n]_q=\frac{q^n-1}{q-1}$ for the $q$-number,
$[n]_{!q}=[n]_q[n-1]_q\cdots[1]_q$ for the $q$-factorial,
$\binom dn_q=\frac{[d]_{!q}}{[n]_{!q}[d-n]_{!q}}$ for the $q$-binomial coefficient, and
\begin{equation}\label{eq:1.1}
S^{\,d}_{n}=\binom{d}{n}_q\,q^{(r-1)n(d-n)}=|\mathcal{X}_n|
\end{equation}
for the number of $n$-spaces in $\mathcal{O}_r^{\,d}$ \cite[Prop. 2.2]{LP25}.

\begin{definition*}[signatures]
A \textbf{signature} is a partition
$\lambda=(\lambda_1\ge\dots\ge\lambda_s>0)$ with $\lambda_1\le r$ and
$s\le d/2$ parts, including the empty partition $\varnothing$. Its
\textbf{size} is $|\lambda|=\sum_i\lambda_i$. We put $s=s(\lambda)$, the number of parts,
for $\lambda\ne\varnothing$, and $s(\varnothing)=1$. The signature $\lambda$
\textbf{lives on the levels} $n=s,s+1,\dots,d-s$, which we call its \textbf{band}; thus
$s(\lambda)$ is the lowest level of the band, and the band of $\varnothing$ consists of
all the levels. So $s=1$ exactly for $\varnothing$ and for the signatures $(t)$ with one part. For $\lambda\ne\varnothing$ the \textbf{exponent} of $\lambda$ is
$\sigma(\lambda)\coloneqq 2rs-|\lambda|$.

\end{definition*}

\subsection{Main results}\label{sec:1.2}

The first result reduces the spectrum to one small matrix per signature.

\begin{restate}{Theorem}{thm:4.4}[block decomposition]
For each signature $\lambda$ there is an
irreducible representation $\chi_\lambda$ of $G=GL_d(\mathcal{O}_r)$, occurring exactly once in
$\ell^2(\mathcal{X}_n)$ for every level $n$ of the band of $\lambda$ and not at all elsewhere,
such that
\[
\operatorname{Spec}(A)\;=\;\biguplus_{\lambda}\;\dim(\chi_\lambda)\cdot
\operatorname{Spec}\bigl(B^{(\lambda)}\bigr),
\]
where $B^{(\lambda)}$ is the symmetric matrix indexed by the band of $\lambda$, with
zero diagonal and off-diagonal entries
\[
B^{(\lambda)}_{mn}\;=\;S^{\,n-s}_{\,m-s}\,
\sqrt{\frac{w^{(\lambda)}_n}{w^{(\lambda)}_m}}\qquad(m<n),
\]
with $s=s(\lambda)$, and $w^{(\lambda)}_n>0$ is the scalar by which $M_{sn}M_{sn}^{\mathsf T}$ acts on
the $\chi_\lambda$-isotypic part of $\ell^2(\mathcal{X}_{s})$; here $M_{mn}$ is the inclusion
matrix between levels $m$ and $n$. When $d$ is even, the signatures with $s=d/2$ live
on the single level $d/2$, and contribute the kernel, of dimension
$|\mathcal{X}_{d/2}|-|\mathcal{X}_{d/2-1}|$.
\end{restate}

The theorem relies on finding the entries $w^{(\lambda)}_n$ and the multiplicities
$\dim\chi_\lambda$, which we compute.

\begin{restate}{Theorem}{thm:6.15}[the entries of the blocks]
Let $\lambda$ be a signature, $s=s(\lambda)$, and let $n$
be a level on its band. Then
\[
w^{(\varnothing)}_n=S^{\,d-1}_{\,n-1}\,S^{\,n}_{1},
\qquad\qquad
w^{(\lambda)}_n\;=\;S^{\,d-2s}_{\,n-s}\;q^{\,(n-s)\,\sigma(\lambda)}\quad(\lambda\ne\varnothing).
\]
In particular $B^{(\varnothing)}$ is the quotient matrix of the level partition, with
entries $\sqrt{S^{\,d-m}_{\,n-m}S^{\,n}_{\,m}}$; its Perron root is the spectral
radius of $A$, which is a simple eigenvalue, and $\mathbb{P}^{d-1}_{fr}(\mathcal{O}_r)$ is connected.
For a signature $(t)$ with one part, $w^{(t)}_n=S^{\,d-2}_{\,n-1}q^{\,(n-1)(2r-t)}$.
\end{restate}

The entries are calculated in \cite{BO12} for $s,n\le d/2$ as eigenvalues of the cellular basis of the Hecke algebra of $\ell^2(\mathcal{X}_n)$, and for $s,n>d/2$ by duality.
For $s<d/2$ and $n>d/2$ the entries are not given by \cite{BO12} and they need more work.
We compute them in two steps: first we reduce the problem to that of one-level eigenvalues, and we then show that up-down operators can be displayed as a sum of down-up operators.

The multiplicities are the dimensions of the
representations $\chi_\lambda$, which were computed by Onn and Stokman \cite{OS06}; in the
labelling of \cite{BO12}, which is the one used here (Section~\ref{sec:3.3}), they read as follows.

\begin{restate}{Theorem}{thm:5.1}[{multiplicities; \cite{OS06}}]
For a signature $(t)$ with one part,
\[
\dim\chi_{(t)}=|\mathcal{X}_1(\mathcal{O}_t)|-|\mathcal{X}_1(\mathcal{O}_{t-1})|
=\begin{cases}q\,[d-1]_q, & t=1,\\[2pt]
[d]_q\,q^{(d-1)(t-2)}\bigl(q^{\,d-1}-1\bigr), & 2\le t\le r.\end{cases}
\]
For $\lambda\ne\varnothing$ in general, let $\lambda'$ be the conjugate partition, and put
$t=1/q$ and $(a;t)_j=\prod_{i=0}^{j-1}(1-t^{i}a)$. Then
\[
\dim\chi_\lambda=t^{\,-(d+1)|\lambda|+2\sum_ii\lambda_i}\,
\frac{\bigl(t^{\,d-\lambda'_1-\lambda'_2+2};t\bigr)_{\lambda'_1+\lambda'_2}\,\bigl(1-t^{\,d-2\lambda'_1+1}\bigr)}
{\prod_{j\ge1}(t;t)_{\lambda'_j-\lambda'_{j+1}}\,\bigl(1-t^{\,d+1}\bigr)} .
\]
In particular $\dim\chi_\lambda$ does not depend on $r$.
\end{restate}

Finally, the blocks determine the spectrum of every subgraph induced on a set of
levels.

\begin{restate}{Theorem}{thm:7.1}[level-set subgraphs]
For $\varnothing\ne S\subseteq\{1,\dots,d-1\}$,
the adjacency matrix $A_S$ of $\mathcal{X}_S$ satisfies
\[
\operatorname{Spec}(A_S)\;=\;\biguplus_{\lambda}\;\dim(\chi_\lambda)\cdot
\operatorname{Spec}\Bigl(B^{(\lambda)}\big|_{\,S\cap\,\mathrm{band}(\lambda)}\Bigr),
\]
the union over the signatures whose band meets $S$, where $B|_T$ denotes the
principal submatrix on the rows and columns indexed by $T$.
\end{restate}

For two levels this is completely explicit.

\begin{restate}{Corollary}{thm:7.2}[two-level subgraphs]
For $1\le m<n\le d-1$,
\[
\operatorname{Spec}\bigl(\mathcal{X}_{m,n}\bigr)
=\biguplus_{\lambda}\ \dim\chi_\lambda\cdot\bigl\{\pm B^{(\lambda)}_{mn}\bigr\}
\ \uplus\ \{0\}^{\,|\mathcal{X}_m|+|\mathcal{X}_n|-2\sum_\lambda\dim\chi_\lambda},
\]
the unions over the signatures with $s(\lambda)\le m$ and $n\le d-s(\lambda)$, where
\[
B^{(\varnothing)}_{mn}=\sqrt{S^{\,d-m}_{\,n-m}\,S^{\,n}_{\,m}},\qquad
B^{(\lambda)}_{mn}\;=\;S^{\,n-s}_{\,m-s}\,
\sqrt{\frac{S^{\,d-2s}_{\,n-s}}{S^{\,d-2s}_{\,m-s}}}\;\;
q^{\,(n-m)\,\sigma(\lambda)/2}\quad(\lambda\ne\varnothing,\ s=s(\lambda)).
\]
\end{restate}

All of the answers above are expressed through $q$, $r$, $d$ and $\lambda$ alone, so
they do not see the ring beyond its residue order and its length (Remark~\ref{thm:7.5}).
Combined with the rigidity theorems of \cite{LP25} this gives:

\begin{restate}{Theorem}{thm:7.6}[isospectral non-isomorphic geometries]
Let $d\ge3$.
\begin{enumerate}
\item If $\mathcal{O}_r=\mathcal{O}/\mathfrak{p}^r$ and
  $\mathcal{O}'_r=\mathcal{O}'/\mathfrak{p}'^{\,r}$ are as above, with the same
  residue order $q$ and the same $r$, then for every nonempty set $S$ of levels the
  graphs $\mathcal{X}_S(\mathcal{O}_r)$ and $\mathcal{X}_S(\mathcal{O}'_r)$ are
  isospectral.
\item Let $p$ be a prime and $r\ge2$. For every set $S$ of at least two levels the graphs
  \[
  \mathcal{X}_S\bigl(\mathbb{Z}/p^r\bigr)\qquad\text{and}\qquad
  \mathcal{X}_S\bigl(\mathbb{F}_p[t]/(t^r)\bigr)
  \]
  are isospectral and non-isomorphic. In particular this holds for the two-level graphs
  $\mathbb{P}^{\,d-1}_{m,n}$, $1\le m<n\le d-1$, and for the full geometries
  $\mathbb{P}^{\,d-1}_{fr}(\mathbb{Z}/p^r)$ and $\mathbb{P}^{\,d-1}_{fr}(\mathbb{F}_p[t]/(t^r))$.
\end{enumerate}
\end{restate}

\subsection{Relation to previous work}\label{sec:1.3}

At $r=1$ the block structure of Theorem~\ref{thm:4.4} is due to Lew \cite{Lew24}, who
block-diagonalizes the weighted upper Laplacian of the flag complex of $\mathbb{F}_q^d$
into blocks of size $d-2m+1$ on the levels $m,\dots,d-m$ with multiplicity
$\binom dm_q-\binom d{m-1}_q$ - the same inventory, with $m=s$ - using Kantor's
inclusion matrices \cite{Kan72}; in substance this is the classical harmonic decomposition
of the subspace lattice. The blocks there are not symmetric, and their eigenvalues are
determined only as $q\to\infty$, where each block tends to a path on its levels. Here
the blocks are symmetric, which is a small improvement. Further advances of Lew's work is in \cite{Vol25}.

The representations $\ell^2(\mathcal{X}_n)$ were studied by Bader and Onn \cite{BO07,BO12}: for
$2n\le d$ they classify the $G$-orbits on $\mathcal{X}_n\times \mathcal{X}_n$, show that $\ell^2(\mathcal{X}_n)$ is
multiplicity-free, label its constituents by partitions, and compute the eigenvalues of
the corresponding Hecke algebra - the character table of the association scheme of
$n$-spaces - in closed form \cite{BO12}. The dimensions of the constituents were computed
by Onn and Stokman \cite{OS06}. We use these results as they stand: the eigenvalues of
\cite{BO12} give the entries of the blocks (Section~\ref{sec:6.1}), and the dimensions of \cite{OS06} give
their multiplicities (Theorem~\ref{thm:5.1}). What we add is the block structure of the whole
multi-level geometry and the identification of the constituents across levels
(Sections~\ref{sec:3.3} and \ref{sec:6.1}); and an additional arguement to move between the lower half and the upper half of the geometry.

\subsection{Organization}\label{sec:1.5}

Section~\ref{sec:2} fixes notation. Section~\ref{sec:3} classifies pairs of summands, and uses the
classification to label the constituents of all the levels by signatures. Section~\ref{sec:4}
proves Theorem~\ref{thm:4.4}, assuming the positivity of the numbers $w^{(\lambda)}_n$.
Section~\ref{sec:dims} gives the dimensions of the constituents, which are the
multiplicities of the blocks (Theorem~\ref{thm:5.1}). Section~\ref{sec:w} computes the numbers
$w^{(\lambda)}_n$: on the lower half of every band from the cellular basis of \cite{BO12}, and
on the whole band by duality; this proves Theorem~\ref{thm:6.15}, completes the proof of
Theorem~\ref{thm:4.4}, and gives Corollary~\ref{thm:6.18}. Section~\ref{sec:7} proves
Theorems~\ref{thm:7.1} and \ref{thm:7.6} and Corollary~\ref{thm:7.2}.

\subsection{Acknowledgements}\label{sec:ack}

The author was supported by the ERC grant No.~882751 (PI: Alex Lubotzky) and by the ISF
grant No.~2990/21 (PI: Ori Parzanchevski).
The author thanks Ori Parzanchevski and Alex Lubotzky for much helpful advice and guidance.

In writing this article, the author used Anthropic's claude models (mainly Opus and Fable). Especially helpful was the reference to \cite{BO12} in the context of the representation theory of $GL_d(\mathcal{O}_r)$. The author also used claude to check the consistency of notation and soundness of logic.

\section{Preliminaries}\label{sec:2}

In this section we fix the notation used throughout, and collect the standard facts
that are used without proof.

\emph{The ring and the spaces.} Let $\mathcal{O}$ be a Dedekind domain, $\mathfrak{p}$ a nonzero prime ideal of residue order $q=|\mathcal{O}/\mathfrak{p}|$, and $r\ge1$. We write $\mathcal{O}_r=\mathcal{O}/\mathfrak{p}^r$ for the finite quotient of $\mathcal{O}$.
As in \cite{LP25}, we denote the image of $\mathfrak{p}$ in $\mathcal{O}_r$ by $\mathfrak{p}$ as well, and fix a uniformizer $\pi\in\mathfrak{p}\setminus\mathfrak{p}^2$, whose image in $\mathcal{O}_r$ we also denote by $\pi$.
Since $\mathcal{O}$ is a Dedekind domain, $\mathcal{O}_r$ is a finite chain ring with maximal ideal
$\mathfrak{p}=\pi\mathcal{O}_r$, group of units
$\mathcal{O}_r^{\times}=\mathcal{O}_r\setminus\mathfrak{p}$, and residue field
$\mathbb{F}_q=\mathcal{O}/\mathfrak{p}$; its ideals form the chain
$\mathcal{O}_r=\mathfrak{p}^0\supsetneq\mathfrak{p}\supsetneq\dots\supsetneq\mathfrak{p}^{r}=0$,
with $\mathfrak{p}^i=\pi^{\,i}\mathcal{O}_r\cong\mathcal{O}_{r-i}$ as modules and $|\mathfrak{p}^{i}|=q^{\,r-i}$. We also have
$\mathcal{O}_r\cong\mathfrak{O}/\mathfrak{m}^r$, where $\mathfrak{O}$ is the completion
of $\mathcal{O}_{\mathfrak{p}}$, the ring of integers of a non-archimedean local field
with residue field $\mathbb{F}_q$ and maximal ideal $\mathfrak{m}$. The latter is the
setting of \cite{BO12} and \cite{OS06}.

We write $\overline{x}$ for reduction modulo $\mathfrak{p}$. Throughout, $d\ge3$ and $V=\mathcal{O}_r^{\,d}$. A
vector $x\in V$ is \textbf{primitive} if it has a coordinate in $\mathcal{O}_r^{\times}$, i.e.
$x\notin\mathfrak{p}V$; equivalently, $x$ belongs to some basis of $V$. More generally, any linearly independent set of vectors in $V$ can be extended to a basis;
a free submodule $M\le V$ of rank $n$, which we will call an $n$-space, is a direct summand of $V$.

\emph{Counts.} We write $[n]_q=\frac{q^n-1}{q-1}$ for the $q$-number,
$[n]_{!q}=[n]_q[n-1]_q\cdots[1]_q$ for the $q$-factorial, and
$\binom kn_q=\frac{[k]_{!q}}{[n]_{!q}[k-n]_{!q}}$ for the $q$-binomial coefficient. The
number of $n$-spaces in $\mathcal{O}_r^{\,k}$ is \cite[Prop. 2.2]{LP25}
\begin{equation}\label{eq:S}
S^{\,k}_{n}=\binom{k}{n}_q\,q^{(r-1)n(k-n)} .
\end{equation}

We use constantly that summands pass to quotients: if $M$ is an
$m$-space, then $U\mapsto U/M$ is an inclusion-preserving bijection from the $n$-spaces
containing $M$ onto the $(n-m)$-spaces of $V/M\cong \mathcal{O}_r^{\,d-m}$. Consequently, for
$m<n$, a fixed $m$-space lies in exactly $S^{\,d-m}_{\,n-m}$ $n$-spaces and a fixed $n$-space
contains exactly $S^{\,n}_{m}$ $m$-spaces.

\emph{Projective geometry.} Let $\mathcal{X}_n$ denote the set of all $n$-spaces in $V$, for $1\le n\le d-1$. We call $\mathcal{X}_n$ a \textbf{level}, and, as in \cite{LP25}, we also say that its vertices have \textbf{color} $n$. By \eqref{eq:S}, $|\mathcal{X}_n|=S^{\,d}_{\,n}$.
The \textbf{free projective geometry} $\mathbb{P}^{d-1}_{fr}(\mathcal{O}_r)$ is the graph with vertex set $\mathcal{X}=\bigsqcup_n\mathcal{X}_n$ in which two spaces are adjacent when one strictly contains the other.
We write $\mathcal{X}_S$ for the induced subgraph on a set $S$ of levels, and $\mathcal{X}_{m,n}=\mathcal{X}_{\{m,n\}}$. The group $G=GL_d(\mathcal{O}_r)$ acts on $\mathcal{X}$ preserving inclusion.

\emph{Matrices.} For $m<n$, $M_{mn}\in\{0,1\}^{\mathcal{X}_m\times \mathcal{X}_n}$ is the \textbf{inclusion matrix},
$M_{mn}(U,W)=1$ if and only if $U\subset W$, and we put $M_{nn}=I$. The adjacency
matrix $A$ on $\ell^2(\mathcal{X})=\bigoplus_n\ell^2(\mathcal{X}_n)$ is the block matrix with $M_{mn}$
above the diagonal, $M_{mn}^{\mathsf T}$ below it, and $0$ on it. For $m\le n$ we write
\[
T^{(m)}_n:=M_{mn}M_{mn}^{\mathsf T},\qquad D^{(n)}_m:=M_{mn}^{\mathsf T}M_{mn},
\]
so that $T^{(m)}_n(U,U')$, for $U,U'\in\mathcal{X}_m$, is the number of $n$-spaces containing
$U+U'$, and $D^{(n)}_m(W,W')$, for $W,W'\in\mathcal{X}_n$, the number of $m$-spaces contained
in $W\cap W'$. The superscript is
the level on which the operator acts. By the quotient bijection,
\begin{equation}\label{eq:2.2}
M_{jm}M_{mn}=S^{\,n-j}_{\,m-j}\,M_{jn}\qquad(j\le m\le n),
\end{equation}
since between a $j$-space and an $n$-space containing it there are exactly
$S^{\,n-j}_{\,m-j}$ $m$-spaces.

\emph{Representations.} For a finite set $Y$, $\ell^2(Y)$ is the space of complex functions
on $Y$ with $\langle f,g\rangle=\sum_yf(y)\overline{g(y)}$; when $G$ acts on $Y$ it is a
unitary \textbf{permutation representation}. The \textbf{$\chi$-isotypic component} of a
representation is the sum of its subrepresentations isomorphic to $\chi$, and the
representation is \textbf{multiplicity-free} if no irreducible occurs in it twice. We
denote by $\mathbbm{1}$ the constant function.

\emph{Schemes.} When $G$ acts on $Y$, the \textbf{orbit matrices} $E_{\mathfrak{a}}$, one for each $G$-orbit
$\mathfrak{a}$ on $Y\times Y$, with $E_{\mathfrak{a}}(y,y')=1$ exactly on the pairs of that orbit, form a basis
of $\operatorname{End}_G\ell^2(Y)$. We call their span the \textbf{scheme algebra} of $Y$,
and the coefficient of $E_{\mathfrak{a}}$ in an operator its \textbf{class value} on $\mathfrak{a}$.

\emph{Partitions and signatures.} A \textbf{partition} of \textbf{length} $m$ is a weakly
decreasing sequence $a=(a_1\ge a_2\ge\dots\ge a_m\ge0)$ of integers with $a_1\le r$; zero
parts are allowed. A \textbf{signature} is a partition
$\lambda=(\lambda_1\ge\dots\ge\lambda_s)$ with no zero parts and $s\le d/2$, including the
empty signature $\varnothing$; removing the zero parts of a partition with at most $d/2$
nonzero parts gives a signature, and a signature with at most $m$ parts, padded with
zeros, is a partition of length $m$. For a signature $\lambda$ we put $s(\lambda)=s$, its
number of parts, if $\lambda\ne\varnothing$, and $s(\varnothing)=1$. Its \textbf{size} is
$|\lambda|=\sum_i\lambda_i$, its \textbf{conjugate} is the vector
$\lambda'=(\lambda'_1,\dots,\lambda'_r)$ with $\lambda'_i=\#\{k:\lambda_k\ge i\}$, and we
put $\lambda'_{r+1}=0$; for $\lambda\ne\varnothing$ its \textbf{exponent} is
$\sigma(\lambda)=2rs-|\lambda|$. The signature $\lambda$ \textbf{lives on the levels}
$n=s(\lambda),\dots,d-s(\lambda)$, which form its \textbf{band}. For signatures
$\lambda,\mu$, $\lambda\subseteq\mu$ means $\lambda_i\le\mu_i$ for all $i$, a missing
part counting as $0$. We write $(c^j)$ for the signature with $j$ parts equal to $c$, and
$(c^j,i)$ for the one with an additional last part $0\le i\le c$, read as $(c^j)$ when
$i=0$.

We use the following four standard facts without proof.

\begin{fact}\label{thm:2.1}
\begin{enumerate}
\item A family of vectors in $\mathcal{O}_r^{\,d}$ whose reductions modulo $\mathfrak{p}$ are
  linearly independent extends to a basis of $\mathcal{O}_r^{\,d}$; so does a linearly
  independent family. In particular primitive vectors extend to bases.
\item For finite $G$-sets $Y,Z$,
  $\dim\operatorname{Hom}_G(\ell^2(Y),\ell^2(Z))$ is the number of $G$-orbits on
  $Y\times Z$; and $\ell^2(Y)$ is multiplicity-free if and only if
  $\operatorname{End}_G\ell^2(Y)$ is commutative.
\item (Schur's lemma) If $\chi$ occurs exactly once in $\ell^2(Z)$, then all $G$-maps
  $\chi\to\ell^2(Z)$ are a scalar multiple of one another.
\item (Smith normal form) Every nonzero element of $\mathcal{O}_r$ is a unit times
  $\pi^{\,i}$ for a unique $i<r$. Consequently, after invertible row and column
  operations, a matrix $\Theta$ over $\mathcal{O}_r$ becomes diagonal with diagonal entries
  $\pi^{\,d_1},\pi^{\,d_2},\dots$, where $0\le d_1\le d_2\le\dots\le r$ and $\pi^{\,r}=0$.
\end{enumerate}
\end{fact}

In a multiplicity-free $\ell^2(Y)$ every isotypic component is irreducible, and by
Fact~\ref{thm:2.1}(3) every operator in the scheme algebra acts on it by a scalar. Conversely, an
operator acting by a scalar on every isotypic component is a combination of the
isotypic projections, which are equivariant; so the scheme algebra of a
multiplicity-free $\ell^2(Y)$ has dimension equal to the number of its constituents.
Finally, since $G$ preserves inclusion, $A$ and all the operators $M_{mn}$, $T^{(m)}_n$,
$D^{(n)}_m$ are $G$-equivariant.

We need one count, a special case of \cite[Prop. 2.3(2)]{LP25}:

\begin{lemma}\label{thm:2.2}
Let $w\in \mathcal{O}_r^{\,d}$ be primitive, $1\le e\le r$ and $0\le n\le d$.
Then the number of $n$-spaces $U\subseteq \mathcal{O}_r^{\,d}$ with $\pi^{\,r-e}w\in U$ is
\[
F(d,n,e)\;=\;S^{\,d-1}_{\,n-1}\;q^{\,(r-e)(d-n)} .
\]
\end{lemma}

\section{Pairs of summands and the constituents}\label{sec:3}

The main geometric input needed for the spectrum is the classification of pairs of
free submodules under $G$. It is due to Bader and Onn \cite[Cl. 2.2, Lem. 2.3]{BO12}, who prove
it using the self-injectivity of $\mathcal{O}_r$; we give an elementary constructive proof, which
also covers summands of different ranks. We then use it to show that every level is
multiplicity-free, and to label the constituents of all the levels by signatures.

\subsection{The normal form}\label{sec:3.1}

\begin{theorem}[normal form for a pair]\label{thm:3.1}
Let $1\le m\le n\le d-1$, $U\in \mathcal{X}_m$ and
$W\in \mathcal{X}_n$. There are a basis of $V$ and integers $r\ge a_1\ge a_2\ge\dots\ge a_m\ge0$
such that, writing $L=\#\{i:a_i\ge1\}$ and listing the basis as
$e_1,\dots,e_m$, $f_1,\dots,f_L$, $g_1,\dots,g_{n-m}$, $h_1,\dots,h_{d-n-L}$,
\[
U=\langle e_1,\dots,e_m\rangle,\qquad
W=\bigl\langle\,e_i+\pi^{\,r-a_i}f_i\ (1\le i\le m)\,\bigr\rangle
+\langle\,g_1,\dots,g_{n-m}\,\rangle ,
\]
where $\pi^{\,r-a_i}f_i$ is read as $0$ when $a_i=0$. In particular $L\le\min(m,d-n)$.
\end{theorem}

\begin{proof}
Write $V=U\oplus C$ with $C$ free of rank $d-m$ (Fact~\ref{thm:2.1}(1)), fix a basis
$u_1,\dots,u_{d-m}$ of $C$, choose a basis $y_1,\dots,y_n$ of $W$, and write (Fact~\ref{thm:2.1}(4))
\[
y_j=x_j+\pi^{\,d_j}u_j,\qquad x_j\in U .
\]

with $0\le d_1\le\dots\le r$, putting
$d_j=r$ for the columns $j>d-m$ (which are zero). If $d_j\ge1$ then $\overline
y_j=\overline x_j\in\overline U$, so at most $m$ of the $d_j$ are positive. Choose a set
$\mathcal G$ of $n-m$ indices with $d_j=0$, let $\mathcal K=\{j_1,\dots,j_m\}$ be its
complement with $d_{j_1}\le\dots\le d_{j_m}$, and put $a_i=r-d_{j_i}$. One can now
replace $x_j$ by $x_j+\pi^{\,d_j}\zeta_j$, independently for each $j$; use it to make
$x_j=0$ for $j\in\mathcal G$, and to choose the $x_j$ with $j\in\mathcal K$ and $d_j=0$
so that $\overline x_{j_1},\dots,\overline x_{j_m}$ is a basis of $\overline U$ (the
ones with $d_j\ge1$ are already independent). Then $x_{j_1},\dots,x_{j_m}$ is a basis
of $U$ by Fact~\ref{thm:2.1}(1) applied in $U$. Setting $e_i=x_{j_i}$, $f_i=u_{j_i}$ (for
$a_i\ge1$), $g_t=u_j$ ($j\in\mathcal G$), and letting the $h$'s be the remaining $u$'s,
we obtain the asserted form. The basis contains $n+L$ vectors among $e,f,g$, so
$L\le d-n$.
\end{proof}

An immediate consequence of the normal form is the following:

\begin{proposition}[the invariant of a pair]\label{thm:3.2}
In the setting and normal form of Theorem~\ref{thm:3.1}, put
$y_i=e_i+\pi^{\,r-a_i}f_i$. Then
\[
U\cap W=\bigoplus_{i=1}^{m}\mathfrak{p}^{a_i}y_i\;\cong\;\prod_{i=1}^{m}\mathcal{O}_{\,r-a_i},
\qquad
U+W\;\cong\;\mathcal{O}_r^{\,n}\times\prod_{a_i\ge1}\mathfrak{p}^{r-a_i}\times\{0\}^{\,d-n-L}.
\]
In particular $\operatorname{rank}(U\cap W)=\#\{i:a_i=0\}$, and two pairs in
$\mathcal{X}_m\times \mathcal{X}_n$ lie in the same $G$-orbit if and only if they have the same
partition $a=(a_1,\dots,a_m)$, which we call the \textbf{invariant} $a(U,W)$ of the pair.
Removing its zero parts gives a signature with $L$ parts.
\end{proposition}

An immediate consequence of the normal form is the following Corollary, which counts the $G$-orbits on pairs of summands of different ranks.

\begin{corollary}[orbit count]\label{thm:3.3}
For $1\le m,n\le d-1$, the $G$-orbits on
$\mathcal{X}_m\times \mathcal{X}_n$ are in bijection with the signatures with at most
$j=\min(m,n,d-m,d-n)$ parts. In particular their number is $\binom{r+j}{j}$.
\end{corollary}

\begin{proof}
We may assume $m\le n$, so that $j=\min(m,d-n)$. By Proposition~\ref{thm:3.2} the
orbits correspond to the invariants which occur, partitions $a$ of length $m$, and by
Theorem~\ref{thm:3.1} these satisfy $L\le d-n$. Conversely every such $a$ occurs: the
configuration of Theorem~\ref{thm:3.1} needs $n+L\le d$ basis vectors, and the generators of $U$
and of $W$ have independent reductions, so they span summands by Fact~\ref{thm:2.1}(1). Removing
the zero parts of $a$ gives a signature with at most $j$ parts, from which $a$ is
recovered by padding with zeros to length $m$, and the number of these signatures is
$\sum_{i=0}^{j}\binom{r+i-1}{i}=\binom{r+j}{j}$.
\end{proof}

For $m=n=1$ the invariant of two lines $L,L'$ is a single integer
$e(L,L')\in\{0,\dots,r\}$, with $L\cap L'\cong\mathfrak{p}^{e}$; more generally, for a line
$L$ and an $n$-space $U$, $\{c\in \mathcal{O}_r:cL\subseteq U\}=\mathfrak{p}^{e}$ for a unique
$e=e(L,U)$, which classifies the pair (Proposition~\ref{thm:3.2}).

\subsection{Multiplicity-freeness}\label{sec:3.2}

We now move on to consider the representation $\ell^2(\mathcal{X}_n)$, and open with a simple lemma:

\begin{lemma}\label{thm:3.4}
For $U,W\in \mathcal{X}_n$ there is an involution $\tau\in G$ with
$\tau U=W$ and $\tau W=U$.
\end{lemma}

\begin{proof}
In the normal form of Theorem~\ref{thm:3.1} with $m=n$, define $\tau(e_i)=e_i+
\pi^{\,r-a_i}f_i$ and $\tau(f_i)=-f_i$ for $i\le L$, and $\tau=\mathrm{id}$ on the other
basis vectors. On each plane $\langle e_i,f_i\rangle$ its matrix is
$\bigl(\begin{smallmatrix}1&0\\ \pi^{\,r-a_i}&-1\end{smallmatrix}\bigr)$, so $\tau\in G$
and $\tau^2=\mathrm{id}$, and $\tau$ exchanges $e_i$ and $y_i$.
\end{proof}

\begin{theorem}[multiplicity-freeness]\label{thm:3.5}
For every $1\le n\le d-1$,
$\operatorname{End}_G\ell^2(\mathcal{X}_n)$ is commutative, so $\ell^2(\mathcal{X}_n)$ is
multiplicity-free, with $\binom{r+j}{j}$ constituents, where $j=\min(n,d-n)$.
\end{theorem}

\begin{proof}
By Lemma~\ref{thm:3.4} every orbit on $\mathcal{X}_n\times \mathcal{X}_n$ is symmetric, so
$E_{\mathfrak{a}}^{\mathsf T}=E_{\mathfrak{a}}$ for every orbit matrix, and every element of the scheme algebra
is symmetric. For $T,T'$ in the algebra, $TT'$ is again in it, so
$TT'=(TT')^{\mathsf T}=T'T$. The rest follows from Fact~\ref{thm:2.1}(2) and Corollary~\ref{thm:3.3}.
\end{proof}

\subsection{One family of constituents for all levels}\label{sec:3.3}

Theorem~\ref{thm:3.5} decomposes each level separately. To assemble the levels we need to know
that their constituents are the same objects where they overlap, and that they overlap
along bands. Write $b(n)=\min(n,d-n)$ and $\mathcal N(j)=\binom{r+j}{j}$.

\begin{theorem}[labelling]\label{thm:3.6}
Let $\mathcal C_n$ be the set of irreducible
representations occurring in $\ell^2(\mathcal{X}_n)$. If $b(m)\le b(n)$ then
$\mathcal C_m\subseteq\mathcal C_n$. Consequently there is a bijection
$\lambda\mapsto\chi_\lambda$ from the signatures onto $\bigcup_n\mathcal C_n$ such that
for every level $n$,
\[
\chi_\lambda\ \text{occurs in}\ \ell^2(\mathcal{X}_n)\iff s(\lambda)\le b(n)
\iff n\ \text{lies on the band of}\ \lambda,
\]
and then it occurs exactly once.
\end{theorem}

\begin{proof}
By Fact~\ref{thm:2.1}(2) and Corollary~\ref{thm:3.3},
$\dim\operatorname{Hom}_G(\ell^2(\mathcal{X}_m),\ell^2(\mathcal{X}_n))=\mathcal N(\min(b(m),b(n)))$, and
since both representations are multiplicity-free this is
$|\mathcal C_m\cap\mathcal C_n|$, while $|\mathcal C_n|=\mathcal N(b(n))$ by
Theorem~\ref{thm:3.5}. If $b(m)\le b(n)$ then $|\mathcal C_m\cap\mathcal C_n|=|\mathcal C_m|$, so
$\mathcal C_m\subseteq\mathcal C_n$. Hence $\mathcal C_n$ depends only on $b(n)$. Putting
also $\mathcal C_0=\mathcal C_d=\{\mathbf 1\}$, the constituents of $\ell^2$ of the
one-point sets $\{0\}$ and $\{V\}$, where $\mathbf 1$ is the trivial representation, we
get the chain
\[
\{\mathbf 1\}=\mathcal C_0\subseteq\mathcal C_1\subseteq\dots\subseteq
\mathcal C_{\lfloor d/2\rfloor}=\mathcal C_{\lceil d/2\rceil}\supseteq\dots\supseteq
\mathcal C_{d-1}\supseteq\mathcal C_d=\{\mathbf 1\},
\]
with $\mathcal C_n=\mathcal C_{d-n}$ and $|\mathcal C_n|=\mathcal N(b(n))$. For
$1\le j\le d/2$, exactly $\mathcal N(j)-\mathcal N(j-1)=\binom{r+j-1}{j}$ irreducibles
appear in passing from $\mathcal C_{j-1}$ to $\mathcal C_j$, which is the number of
signatures with exactly $j$ parts; any bijection between these layers has the required
property.
\end{proof}

Theorem~\ref{thm:3.6} leaves the bijection free within each layer, and Section~\ref{sec:4} needs nothing
more, since Theorem~\ref{thm:4.4} is a statement about all signatures at once. To compute the
blocks we fix the labelling of Bader and Onn \cite[\S2]{BO12}: for $2m\le d$ they index the
constituents of $\ell^2(\mathcal{X}_m)$ by the signatures $\lambda\subseteq(r^m)$, i.e. by those
with at most $m$ parts, $\lambda\mapsto\beta_m(\lambda)$
(Theorem~\ref{thm:6.2}). By Lemma~\ref{thm:6.3}, $\beta_m(\lambda)$ is the same
representation at every level $m$ with $s(\lambda)\le m\le d/2$, and it is new at level
$s(\lambda)$. We put $\chi_\lambda:=\beta_m(\lambda)$; by Theorem~\ref{thm:3.6} this is a
labelling as there, and it determines $\chi_\lambda$ at the levels above $d/2$ as well.

\begin{corollary}\label{thm:3.7}
(1) For $m$ outside the band of $\lambda$, the
$\chi_\lambda$-isotypic part of $\ell^2(\mathcal{X}_m)$ is $0$. (2) For every level $n$,
\begin{equation}\label{eq:3.1}
\sum_{\lambda:\ s(\lambda)\le b(n)}\dim\chi_\lambda\;=\;|\mathcal{X}_n| .
\end{equation}
\end{corollary}

\begin{proof}
Both assertions are immediate from Theorems~\ref{thm:3.5} and \ref{thm:3.6}, (2) being the
dimension count of the decomposition of $\ell^2(\mathcal{X}_n)$ into its constituents.
\end{proof}

For example, for $d=6$ the levels are $1,\dots,5$, and the signatures live as follows:

\begin{center}
\begin{tabular}{ccccccc}
signature $\backslash$ level & 1 & 2 & 3 & 4 & 5 & block size \\
\hline
$\varnothing$ & $\bullet$ & $\bullet$ & $\bullet$ & $\bullet$ & $\bullet$ & $5$ \\
$(t)$, $1\le t\le r$ & $\bullet$ & $\bullet$ & $\bullet$ & $\bullet$ & $\bullet$ & $5$ \\
$(t,u)$ &  & $\bullet$ & $\bullet$ & $\bullet$ &  & $3$ \\
$(t,u,v)$ &  &  & $\bullet$ &  &  & $1$ (kernel) \\
\end{tabular}
\end{center}

\section{The block decomposition}\label{sec:4}

In this section we prove Theorem~\ref{thm:4.4}, assuming that the numbers $w^{(\lambda)}_n$ are
positive; this is established in Section~\ref{sec:w}. Throughout, $\lambda$ is a signature, $s$
the lowest level of its band, and $W^{(n)}_\lambda$ the $\chi_\lambda$-isotypic part of
$\ell^2(\mathcal{X}_n)$, which is irreducible for $n$ on the band (Theorems~\ref{thm:3.5} and \ref{thm:3.6}). We
write $d_\lambda=\dim\chi_\lambda$.

\begin{lemma}\label{thm:4.1}
For every signature $\lambda$ and every level $n\ge s$, the operator
$T^{(s)}_n$ acts on $W^{(s)}_\lambda$ by a scalar, which we denote by
$w^{(\lambda)}_n$. Moreover $w^{(\lambda)}_{s}=1$, and $w^{(\lambda)}_n=0$ for
$n>d-s$.
\end{lemma}

\begin{proof}
$T^{(s)}_n$ is equivariant and $W^{(s)}_\lambda$ is irreducible, so this
is Fact~\ref{thm:2.1}(3). For $n=s$ we have $T^{(s)}_{s}=I$. For $n>d-s$ the
representation $\chi_\lambda$ does not occur in $\ell^2(\mathcal{X}_n)$, so
$M_{sn}^{\mathsf T}$ kills $W^{(s)}_\lambda$.
\end{proof}

We say that $\lambda$ satisfies \textbf{condition (P)} if $w^{(\lambda)}_n>0$ for every
level $n$ on its band.

\begin{proposition}\label{thm:4.2}
Let $n$ be a level on the band of $\lambda$. If $\lambda$
satisfies (P), then $M_{sn}^{\mathsf T}$ is injective on $W^{(s)}_\lambda$.
\end{proposition}

\begin{proof}
For $\varphi\in W^{(s)}_\lambda$,
$\lVert M_{sn}^{\mathsf T}\varphi\rVert^2=\langle T^{(s)}_n\varphi,\varphi\rangle=
w^{(\lambda)}_n\lVert\varphi\rVert^2$.
\end{proof}

\begin{proposition}[the block of a signature]\label{thm:4.3}
Assume (P) for $\lambda$, let
$\varphi\in W^{(s)}_\lambda$ be a unit vector, and put
$\psi_n=M_{sn}^{\mathsf T}\varphi/\sqrt{w^{(\lambda)}_n}$ for $n$ on the band. Then
the $\psi_n$ are orthonormal, their span is $A$-invariant, and in the basis
$(\psi_n)_n$ the operator $A$ acts by the matrix $B^{(\lambda)}$ of Theorem~\ref{thm:4.4}:
\[
\langle A\psi_n,\psi_m\rangle
=\beta_{mn}\sqrt{\frac{w^{(\lambda)}_n}{w^{(\lambda)}_m}},
\qquad
\beta_{mn}=S^{\,n-s}_{\,m-s}\qquad(m<n).
\]
Moreover $A\psi_n$ has no component on the levels outside the band.
\end{proposition}

\begin{proof}
The $\psi_n$ are unit vectors by the proof of Proposition~\ref{thm:4.2}, and lie in
different levels. The component of $A\psi_n$ on level $m$ is $M_{mn}\psi_n$ for $m<n$
and $M_{nm}^{\mathsf T}\psi_n$ for $m>n$; by equivariance it lies in
$W^{(m)}_\lambda$, which is $0$ if $m$ is off the band (Corollary~\ref{thm:3.7}). Let $m$ be on
the band. For $m>n$, \eqref{eq:2.2} gives
$M_{nm}^{\mathsf T}M_{sn}^{\mathsf T}=\beta_{nm}M_{sm}^{\mathsf T}$, so the
component is a multiple of $\psi_m$. For $m<n$, the map
$\varphi'\mapsto M_{mn}M_{sn}^{\mathsf T}\varphi'$ is a $G$-map from
$W^{(s)}_\lambda$ to $\ell^2(\mathcal{X}_m)$, hence by Fact~\ref{thm:2.1}(3) a multiple of
$M_{sm}^{\mathsf T}$, and again $M_{mn}\psi_n\in\mathbb{C}\psi_m$. The coefficient is
\[
\langle M_{mn}\psi_n,\psi_m\rangle
=\frac{\langle M_{sm}M_{mn}M_{sn}^{\mathsf T}\varphi,\varphi\rangle}
{\sqrt{w^{(\lambda)}_mw^{(\lambda)}_n}}
=\beta_{mn}\,\frac{\langle T^{(s)}_n\varphi,\varphi\rangle}
{\sqrt{w^{(\lambda)}_mw^{(\lambda)}_n}}
=\beta_{mn}\sqrt{\frac{w^{(\lambda)}_n}{w^{(\lambda)}_m}},
\]
using \eqref{eq:2.2} once more. The diagonal vanishes since the levels carry no edges.
\end{proof}

Note that for $m<n$ the operator $M_{sm}M_{mn}M_{sn}^{\mathsf T}$ is not a scalar on
$\ell^2(\mathcal{X}_{s})$; it becomes one only after restriction to the irreducible
$W^{(s)}_\lambda$, and this is the only point at which the invariance of the span of
the $\psi_n$ uses more than counting.

\begin{theorem}[block decomposition]\label{thm:4.4}
For each signature $\lambda$ there is an
irreducible representation $\chi_\lambda$ of $G$, occurring exactly once in
$\ell^2(\mathcal{X}_n)$ for every level $n$ of the band of $\lambda$ and not at all elsewhere,
such that
\[
\operatorname{Spec}(A)\;=\;\biguplus_{\lambda}\;\dim(\chi_\lambda)\cdot
\operatorname{Spec}\bigl(B^{(\lambda)}\bigr),
\]
where $B^{(\lambda)}$ is the symmetric matrix indexed by the band of $\lambda$, with
zero diagonal and off-diagonal entries
\[
B^{(\lambda)}_{mn}\;=\;S^{\,n-s}_{\,m-s}\,
\sqrt{\frac{w^{(\lambda)}_n}{w^{(\lambda)}_m}}\qquad(m<n),
\]
with $s=s(\lambda)$, and $w^{(\lambda)}_n>0$ is the scalar of Lemma~\ref{thm:4.1}. When $d$ is even, the signatures with $s=d/2$ live
on the single level $d/2$, and contribute the kernel, of dimension
$|\mathcal{X}_{d/2}|-|\mathcal{X}_{d/2-1}|$.
\end{theorem}

\begin{proof}[Proof of Theorem~\ref{thm:4.4}, assuming (P) for every signature]
By Theorems~\ref{thm:3.5} and \ref{thm:3.6},
$\ell^2(\mathcal{X}_n)=\bigoplus_{\lambda}W^{(n)}_\lambda$, the sum over the $\lambda$ whose band
contains $n$, with $\dim W^{(n)}_\lambda=d_\lambda$; this is \eqref{eq:3.1}. Hence
$\ell^2(\mathcal{X})=\bigoplus_\lambda Y_\lambda$, where
$Y_\lambda=\bigoplus_nW^{(n)}_\lambda$ is the $\chi_\lambda$-isotypic component of
$\ell^2(\mathcal{X})$, and $AY_\lambda\subseteq Y_\lambda$ by equivariance. Fix an orthonormal
basis $\varphi_1,\dots,\varphi_{d_\lambda}$ of $W^{(s)}_\lambda$, and let $U_i$ be the
span of the vectors $\psi^i_n$ obtained from $\varphi_i$ as in Proposition~\ref{thm:4.3}. Each
$U_i$ is $A$-invariant, and $A|_{U_i}$ is $B^{(\lambda)}$. The $U_i$ are mutually
orthogonal, since
$\langle M_{sn}^{\mathsf T}\varphi_i,M_{sn}^{\mathsf T}\varphi_j\rangle=
w^{(\lambda)}_n\langle\varphi_i,\varphi_j\rangle=0$ for $i\ne j$, and their dimensions
add up to $\dim Y_\lambda$. Thus $A$ is the orthogonal sum, over all signatures, of
$d_\lambda$ copies of $B^{(\lambda)}$. If $d$ is even, the signatures with $s=d/2$ live
only on the level $d/2$, so their blocks are $1\times1$ zero matrices, and by \eqref{eq:3.1} at
$n=d/2$ and $n=d/2-1$ their total multiplicity is $|\mathcal{X}_{d/2}|-|\mathcal{X}_{d/2-1}|$.
\end{proof}

The proof uses no value of $w^{(\lambda)}_n$, only its sign; in particular Theorem~\ref{thm:4.4}
determines the size and multiplicity of every block independently of its entries.

\section{Dimensions}\label{sec:dims}

By Theorem~\ref{thm:4.4}, the block $B^{(\lambda)}$ occurs in the spectrum of $A$ with
multiplicity $\dim\chi_\lambda$. These dimensions were computed by Onn and Stokman \cite{OS06}.

\begin{theorem}[{multiplicities; \cite{OS06}}]\label{thm:5.1}
For a signature $(t)$ with one part,
\[
\dim\chi_{(t)}=|\mathcal{X}_1(\mathcal{O}_t)|-|\mathcal{X}_1(\mathcal{O}_{t-1})|
=\begin{cases}q\,[d-1]_q, & t=1,\\[2pt]
[d]_q\,q^{(d-1)(t-2)}\bigl(q^{\,d-1}-1\bigr), & 2\le t\le r.\end{cases}
\]
For $\lambda\ne\varnothing$ in general, put $t=1/q$, $(a;t)_j=\prod_{i=0}^{j-1}(1-t^{i}a)$ and $(t;t)_n=\prod_{i=1}^{n}(1-t^{i})$. Then
\[
\dim\chi_\lambda=t^{\,-(d+1)|\lambda|+2\sum_ii\lambda_i}\,
\frac{\bigl(t^{\,d-\lambda'_1-\lambda'_2+2};t\bigr)_{\lambda'_1+\lambda'_2}\,\bigl(1-t^{\,d-2\lambda'_1+1}\bigr)}
{\prod_{j\ge1}(t;t)_{\lambda'_j-\lambda'_{j+1}}\,\bigl(1-t^{\,d+1}\bigr)} .
\]
In particular $\dim\chi_\lambda$ does not depend on $r$.
\end{theorem}

\begin{proof}[Proof of Theorem~\ref{thm:5.1}]
Onn and Stokman \cite[Thm. 4.5, (4.13)]{OS06} give the dimension of the constituent
$U_\lambda$ of the Grassmann representation, in the labelling of \cite{BO12} used in
Section~\ref{sec:3.3}: if $\lambda$ has $s$ parts and $\rho=(s-1,\dots,1,0)$, then
\[
\dim\chi_\lambda=t^{\,-(d-2s+1)|\lambda|-2\langle\rho,\lambda\rangle}\,
\frac{(t;t)_s\,\bigl(t^{\,d-\lambda'_1-\lambda'_2+2};t\bigr)_{\lambda'_1+\lambda'_2}\,\bigl(1-t^{\,d-2\lambda'_1+1}\bigr)}
{(t;t)_{s-\lambda'_1}\,\bigl(t^{\,s-\lambda'_1+1};t\bigr)_{\lambda'_1}\,\prod_{j\ge1}(t;t)_{\lambda'_j-\lambda'_{j+1}}\,\bigl(1-t^{\,d+1}\bigr)} .
\]
Since $\lambda'_1=s$, we have $(t;t)_{s-\lambda'_1}=1$, $\bigl(t^{\,s-\lambda'_1+1};t\bigr)_{\lambda'_1}=(t;t)_s$
and $\langle\rho,\lambda\rangle=s|\lambda|-\sum_ii\lambda_i$, which gives the formula stated. It contains no $r$.
\end{proof}

\section{\texorpdfstring{The numbers $w^{(\lambda)}_n$}{The numbers w(lambda,n)}}\label{sec:w}

By Lemma~\ref{thm:4.1}, $w^{(\lambda)}_n$ is the eigenvalue of
$T^{(s)}_n=M_{sn}M_{sn}^{\mathsf T}$ on the $\chi_\lambda$-isotypic part of
$\ell^2(\mathcal{X}_{s})$, at the bottom of the band. The same number is the eigenvalue of
$D^{(n)}_{s}=M_{sn}^{\mathsf T}M_{sn}$ on the $\chi_\lambda$-isotypic part of
$\ell^2(\mathcal{X}_n)$. This holds for any pair of levels, not only for the bottom one:

\begin{lemma}\label{thm:w-down}
Let $\lambda$ be a signature and $j<n$ levels. Then $D^{(n)}_{j}$ acts on
$W^{(n)}_\lambda$ and $T^{(j)}_n$ acts on $W^{(j)}_\lambda$ by the same scalar, and this
scalar is $0$ unless $j$ and $n$ both lie on the band of $\lambda$. In particular, with
$s=s(\lambda)$, $D^{(n)}_{s}$ acts on $W^{(n)}_\lambda$ by $w^{(\lambda)}_n$ for every level
$n$ on the band.
\end{lemma}

\begin{proof}
Let $M=M_{jn}$. Being equivariant, $M$ maps $W^{(n)}_\lambda$ into $W^{(j)}_\lambda$ and
$M^{\mathsf T}$ maps $W^{(j)}_\lambda$ into $W^{(n)}_\lambda$, and
$T^{(j)}_nM=MM^{\mathsf T}M=MD^{(n)}_{j}$. By Fact~\ref{thm:2.1}(3), $D^{(n)}_{j}$ acts on
$W^{(n)}_\lambda$ by a scalar $w$ and $T^{(j)}_n$ acts on $W^{(j)}_\lambda$ by a scalar $t$
(either space may be $0$, and then its scalar is arbitrary). If $M\psi\ne0$ for some
$\psi\in W^{(n)}_\lambda$, then both spaces are nonzero, so $j$ and $n$ lie on the band, and
\[
tM\psi=T^{(j)}_nM\psi=MD^{(n)}_{j}\psi=wM\psi
\]
gives $t=w$. Otherwise $M$ kills $W^{(n)}_\lambda$, and we may take $w=t=0$: for
$\psi\in W^{(n)}_\lambda$ we have
$\langle D^{(n)}_{j}\psi,\psi\rangle=\lVert M\psi\rVert^2=0$, and for
$\varphi\in W^{(j)}_\lambda$ we have $M^{\mathsf T}\varphi\in W^{(n)}_\lambda$, so
$\langle T^{(j)}_n\varphi,\varphi\rangle=\lVert M^{\mathsf T}\varphi\rVert^2
=\langle MM^{\mathsf T}\varphi,\varphi\rangle=0$. If $j$ or $n$ is off the band, one of the
two spaces is $0$, so $M$ kills $W^{(n)}_\lambda$ and we are in this second case. The last
assertion is the case $j=s$, by the definition of $w^{(\lambda)}_n$.
\end{proof}

The operator $D^{(n)}_{s}$, which counts the $s$-spaces in the intersection of
two $n$-spaces, is an element of the cellular basis of $\operatorname{End}_G\ell^2(\mathcal{X}_n)$,
whose eigenvalues Bader and Onn computed in closed form for $2n\le d$ \cite{BO12}.
Section~\ref{sec:6.1} makes the translation from their notation to ours and shows that in our case the formula simplifies.
Sections~\ref{sec:6.2} and \ref{sec:6.3} extend the result to the whole
band by duality; Section~\ref{sec:6.4} assembles the blocks and completes the proof of
Theorem~\ref{thm:4.4}; and Section~\ref{sec:6.5} writes every block as one universal
matrix.

\subsection{The lower half of the band}\label{sec:6.1}

\emph{The cellular basis.} For a signature $\mu$, a submodule of $V$
is of \textbf{type} $\mu$ if it is isomorphic to $\bigoplus_i\mathcal{O}_r/\mathfrak{p}^{\mu_i}$.
For a level $m$ let $M_{\mu m}$ be the inclusion matrix between the submodules of type
$\mu$ and $\mathcal{X}_m$, with $M_{\mu m}(Z,U)=1$ if and only if $Z\subseteq U$, and put
\[
\mathbf c^{(m)}_\mu:=M_{\mu m}^{\mathsf T}M_{\mu m},\qquad\text{so that}\qquad
\mathbf c^{(m)}_\mu(U,U')=\#\{Z\subseteq U\cap U':\ Z\ \text{of type}\ \mu\}.
\]
For $\mu\subseteq(r^m)$ these operators form the \textbf{cellular basis} of
$\operatorname{End}_G\ell^2(\mathcal{X}_m)$ of \cite{BO12}. A submodule of type $(r^j)$ is a free
submodule of rank $j$, that is, a $j$-space (Section~\ref{sec:2}); hence
$M_{(r^j)\,m}=M_{jm}$ and
\begin{equation}\label{eq:6.1}
\mathbf c^{(m)}_{(r^j)}=D^{(m)}_j\qquad(j\le m).
\end{equation}

\begin{theorem}[{\cite[Thm. 2, Cor. 3, Thm. 4, Cor. 6]{BO12}}]\label{thm:6.2}
Let $2m\le d$. There is a bijection $\lambda\mapsto\beta_m(\lambda)$ from the signatures
$\lambda\subseteq(r^m)$ onto the constituents of $\ell^2(\mathcal{X}_m)$ such that for every
$\mu\subseteq(r^m)$ the operator $\mathbf c^{(m)}_\mu$ acts on the
$\beta_m(\lambda)$-isotypic part of $\ell^2(\mathcal{X}_m)$ by a scalar
$A^{(m)}(\mu,\lambda)$, which is $0$ unless $\lambda\subseteq\mu$, and for
$\lambda\subseteq\mu$
\[
A^{(m)}(\mu,\lambda)=q^{\langle\phi'-\mu',\,\Phi'-\phi'+\mu'-\lambda'\rangle}
\genfrac{[}{]}{0pt}{}{\phi'_1-\lambda'_1}{\phi'_1-\mu'_1}
\genfrac{[}{]}{0pt}{}{\Phi'_1-\mu'_1-\lambda'_1}{\phi'_1-\mu'_1}
\prod_{i\ge1}\genfrac{[}{]}{0pt}{}{\mu'_i-\lambda'_{i+1}}{\mu'_i-\mu'_{i+1}}.
\]
Here $\phi=(r^m)$ and $\Phi=(r^d)$ are the types of an $m$-space and of $V$, so that
$\phi'=(m,\dots,m)$ and $\Phi'=(d,\dots,d)$; the inner product is that of vectors of length
$r$; and
\[
\genfrac{[}{]}{0pt}{}{N}{K}=\frac{\prod_{k=1}^{N}(1-q^{-k})}
{\prod_{k=1}^{K}(1-q^{-k})\,\prod_{k=1}^{N-K}(1-q^{-k})}=q^{-K(N-K)}\binom NK_q,
\]
which is $0$ unless $0\le K\le N$.
\end{theorem}

For $\lambda\subseteq\mu\subseteq(r^m)$ every factor $\genfrac{[}{]}{0pt}{}{N}{K}$ of the formula has
$0\le K\le N$, the second because $\lambda'_1\le m\le d-m$, and so is nonzero. Hence
\begin{equation}\label{eq:6.2}
A^{(m)}(\mu,\lambda)\ne0\iff\lambda\subseteq\mu .
\end{equation}

Theorem~\ref{thm:6.2} indexes the constituents level by level. The following lemma
shows that the indices agree across the levels $m\le d/2$.

\begin{lemma}[the indices along the lower half]\label{thm:6.3}
Let $j\le n\le d/2$.
\begin{enumerate}
\item For every signature $\mu\subseteq(r^j)$ there is a constant $c>0$ with
  $M_{\mu j} M_{jn}=c\,M_{\mu n}$.
\item The constituents of $\ell^2(\mathcal{X}_j)$ are the $\beta_n(\lambda)$ with
  $\lambda\subseteq(r^j)$, and $\beta_n(\lambda)=\beta_j(\lambda)$ for these $\lambda$.
\end{enumerate}
\end{lemma}

\begin{proof}
(1) The entry of $M_{\mu j} M_{jn}$ at $(Z,W)$ is the number of $j$-spaces $U$ with
$Z\subseteq U\subseteq W$, and it vanishes unless $Z\subseteq W$. Let $Z\subseteq W$,
with $Z$ of type $\mu$ and $W\in \mathcal{X}_n$. Applying Fact~\ref{thm:2.1}(4) to the matrix which
expresses a generating set of $Z$ in a basis of $W$, we obtain a basis $f_1,\dots,f_n$
of $W$ with $Z=\bigoplus_{i\le\ell}\mathfrak{p}^{r-\mu_i}f_i$, where $\ell\le j$
is the number of parts of $\mu$. Extending it to a basis of $V$
(Fact~\ref{thm:2.1}(1)) shows that $G$ is transitive on such pairs $(Z,W)$, so that the
entry is a constant $c$, and $U=\langle f_1,\dots,f_j\rangle$ shows that $c>0$.

(2) Let $\lambda\subseteq(r^j)$ and $\chi=\beta_n(\lambda)$. By \eqref{eq:6.2},
$\mathbf c^{(n)}_\lambda$ acts on the $\chi$-isotypic part of $\ell^2(\mathcal{X}_n)$ by a
nonzero scalar, so $M_{\lambda n}$ does not vanish there, and by (1) neither does
$M_{jn}$. Since $M_{jn}$ is equivariant, $\chi$ occurs in $\ell^2(\mathcal{X}_j)$, and
$M_{jn}$ maps the irreducible $\chi$-isotypic part of $\ell^2(\mathcal{X}_n)$ isomorphically
onto that of $\ell^2(\mathcal{X}_j)$. This gives $\binom{r+j}{j}$
distinct constituents of $\ell^2(\mathcal{X}_j)$, which by Theorem~\ref{thm:3.5} are all of them;
this proves the first assertion of (2). Now let $\chi=\beta_j(\kappa)$. For
$\mu\subseteq(r^j)$, by (1) and since $M_{jn}$ maps the $\chi$-isotypic part of
$\ell^2(\mathcal{X}_n)$ onto that of $\ell^2(\mathcal{X}_j)$, the operator $M_{\mu n}$
vanishes on the first if and only if $M_{\mu j}$ vanishes on the second; that is,
$\mathbf c^{(n)}_\mu$ and $\mathbf c^{(j)}_\mu$ vanish on $\chi$ together. By \eqref{eq:6.2}, $\lambda\subseteq\mu$ if and only if
$\kappa\subseteq\mu$, for every $\mu\subseteq(r^j)$, and taking $\mu=\lambda$ and
$\mu=\kappa$ gives $\kappa=\lambda$.
\end{proof}

Thus for a signature $\lambda$ the representation $\beta_m(\lambda)$ is the same for all
levels $s(\lambda)\le m\le d/2$, and by (2) it occurs in $\ell^2(\mathcal{X}_{s(\lambda)})$ but not
in $\ell^2(\mathcal{X}_{s(\lambda)-1})$. This is the labelling $\chi_\lambda$ of Section~\ref{sec:3.3},
and we do not use $\beta_m$ again: Theorem~\ref{thm:6.2} says that for $2m\le d$ and
$\lambda,\mu\subseteq(r^m)$, $\mathbf c^{(m)}_\mu$ acts on $W^{(m)}_\lambda$ by
$A^{(m)}(\mu,\lambda)$.

\begin{proposition}[the lower half of the band]\label{thm:6.4}
Let $\lambda$ be a signature, and let $s\le n\le d/2$. Then
$w^{(\lambda)}_n=A^{(n)}\bigl((r^{s}),\lambda\bigr)$; explicitly,
\[
w^{(\varnothing)}_n=S^{\,d-1}_{\,n-1}\,S^{\,n}_{1},
\qquad\qquad
w^{(\lambda)}_n=S^{\,d-2s}_{\,n-s}\;q^{\,(n-s)\,\sigma(\lambda)}\quad(\lambda\ne\varnothing).
\]
In particular $w^{(\lambda)}_n>0$.
\end{proposition}

\begin{proof}
By Lemma~\ref{thm:w-down}, $w^{(\lambda)}_n$ is the scalar by which $D^{(n)}_{s}$ acts on
$W^{(n)}_\lambda$. By \eqref{eq:6.1}, $D^{(n)}_{s}=\mathbf c^{(n)}_{(r^{s})}$, and
$\lambda\subseteq(r^s)$. So by Theorem~\ref{thm:6.2},
$w^{(\lambda)}_n=A^{(n)}((r^{s}),\lambda)$. In the formula take $m=n$ and $\mu=(r^s)$, so
that $\phi'=(n,\dots,n)$, $\mu'=(s,\dots,s)$ and $\Phi'=(d,\dots,d)$. Let
$\lambda\ne\varnothing$; then $\lambda'_1=s$ and $\sum_i\lambda'_i=|\lambda|$. Most of the
terms cancel, and by \eqref{eq:S} the formula simplifies to
\[
w^{(\lambda)}_n=S^{\,d-2s}_{\,n-s}\,q^{\,(n-s)\sigma(\lambda)} .
\]
For $\lambda=\varnothing$ we have $s=1$ and $\lambda'=0$: the exponent is $r(n-1)(d-n+1)$,
the two middle factors are $\genfrac{[}{]}{0pt}{}{n}{n-1}=q^{-(n-1)}[n]_q$ and
$\genfrac{[}{]}{0pt}{}{d-1}{n-1}=q^{-(n-1)(d-n)}\binom{d-1}{n-1}_q$, the product is $1$, and the result is
$q^{\,(r-1)(n-1)(d-n+1)}[n]_q\binom{d-1}{n-1}_q=S^{\,d-1}_{\,n-1}S^{\,n}_1$.
\end{proof}

\begin{remark}\label{thm:6.5}
The answer depends on $\lambda$ only through $s$ and $|\lambda|$: the vector
$\phi'-\mu'$ is the constant $(n-s)$, against which $\lambda'$ contributes only its
size. This is the source of the universality of Corollary~\ref{thm:6.18}.
\end{remark}

\subsection{Self-duality}\label{sec:6.2}

The band of a signature with $s$ parts extends up to the level $d-s$, beyond the range
$2m\le d$ of Theorem~\ref{thm:6.2}. We cover the levels above $d/2$ with the map
$U\mapsto U^\perp$, after showing that it carries every constituent to itself. Let
$\langle\ ,\ \rangle$ be the standard bilinear form on $V$, let
$U^\perp=\{w\in V:\langle u,w\rangle=0\ \text{for all}\ u\in U\}$, and let
$\theta(g)=(g^{\mathsf T})^{-1}$, an automorphism of $G$.

\begin{fact}\label{thm:6.6}
For $U\in \mathcal{X}_m$, $U^\perp\in \mathcal{X}_{d-m}$ and $(U^\perp)^\perp=U$; for
$U\in \mathcal{X}_m$ and $W\in \mathcal{X}_n$, $U\subseteq W$ if and only if $W^\perp\subseteq U^\perp$;
and $(gU)^\perp=\theta(g)\,U^\perp$ for $g\in G$. In particular $U\mapsto U^\perp$ is an
automorphism of the graph $\mathbb{P}^{d-1}_{fr}(\mathcal{O}_r)$ which maps $\mathcal{X}_m$ onto $\mathcal{X}_{d-m}$.
\end{fact}

\begin{lemma}[duality preserves the invariant]\label{thm:6.7}
Let $m\le n$ with $m+n=d$. For $U\in \mathcal{X}_m$ and $W\in \mathcal{X}_n$, the pair
$(W^\perp,U^\perp)\in \mathcal{X}_m\times \mathcal{X}_n$ has the same invariant as $(U,W)$.
\end{lemma}

\begin{proof}
By Proposition~\ref{thm:3.2} the invariant of a pair in $\mathcal{X}_m\times \mathcal{X}_n$ is determined
by the isomorphism type of the intersection, and $W^\perp\cap U^\perp=(U+W)^\perp$; so we
show that $(U+W)^\perp\cong U\cap W$. Let $\varphi\colon U\to V/W$ be the restriction of the
quotient map. Its kernel is $U\cap W$ and its image is $(U+W)/W\cong U/(U\cap W)$ (the second
isomorphism theorem), so its cokernel is $V/(U+W)$. Both $U$ and $V/W$ are free of rank
$m=d-n$, so by Fact~\ref{thm:2.1}(4) $\varphi$ is $\operatorname{diag}(\pi^{\,c_1},\dots,\pi^{\,c_m})$
in suitable bases, and
\[
U\cap W=\ker\varphi\cong\bigoplus_{i=1}^{m}\mathfrak{p}^{\,r-c_i}\cong\bigoplus_{i=1}^{m}\mathcal{O}_r/\mathfrak{p}^{\,c_i}
\cong\operatorname{coker}\varphi=V/(U+W).
\]
Finally, every homomorphism $V\to\mathcal{O}_r$ is $\langle\,\cdot\,,v\rangle$ for a unique $v\in V$,
and it vanishes on $U+W$ if and only if $v\in(U+W)^\perp$. Hence
$(U+W)^\perp\cong\operatorname{Hom}(V/(U+W),\mathcal{O}_r)\cong V/(U+W)$, since
$\operatorname{Hom}(\mathcal{O}_r/\mathfrak{p}^{\,c},\mathcal{O}_r)\cong\mathfrak{p}^{\,r-c}\cong\mathcal{O}_r/\mathfrak{p}^{\,c}$.
\end{proof}

\begin{proposition}[self-duality]\label{thm:6.8}
Let $\Phi$ be the unitary involution of $\ell^2(\mathcal{X})$ given by $(\Phi f)(U)=f(U^\perp)$.
Then $\Phi W^{(j)}_\lambda=W^{(d-j)}_\lambda$ for every signature $\lambda$ and every level $j$.
\end{proposition}

\begin{proof}
By Fact~\ref{thm:6.6}, $\Phi$ maps $\ell^2(\mathcal{X}_j)$ onto $\ell^2(\mathcal{X}_{d-j})$ and
$\Phi(\theta(g)\cdot f)=g\cdot\Phi f$, so it carries $\chi$-isotypic parts to
$\chi\circ\theta$-isotypic parts: $\Phi W^{(j)}_\lambda=W^{(d-j)}_\mu$ with
$\chi_\mu\cong\chi_\lambda\circ\theta$, and we must show $\mu=\lambda$. Let $m=\lfloor d/2\rfloor$ and
$n=d-m$; then $\chi_\lambda$ occurs once in each of $\ell^2(\mathcal{X}_m)$ and $\ell^2(\mathcal{X}_n)$
(Theorem~\ref{thm:3.6}). By Lemma~\ref{thm:6.7}, $\Phi E_{\mathfrak{a}}\Phi=E_{\mathfrak{a}}^{\mathsf T}$
for every orbit matrix on $\mathcal{X}_m\times\mathcal{X}_n$, hence (Fact~\ref{thm:2.1}(2))
$\Phi E\Phi=E^{\mathsf T}$ for every $G$-map $E\colon\ell^2(\mathcal{X}_n)\to\ell^2(\mathcal{X}_m)$.
Take $E$ mapping $W^{(n)}_\lambda$ isomorphically onto $W^{(m)}_\lambda$ and vanishing on the other
isotypic parts. Then $EE^{\mathsf T}=\Phi(E^{\mathsf T}E)\Phi$ is injective on $W^{(m)}_\lambda$ and
vanishes on the orthogonal complement of $\Phi W^{(n)}_\lambda=W^{(m)}_\mu$; hence $\mu=\lambda$.
\end{proof}

By the first step of the proof, Proposition~\ref{thm:6.8} says that every constituent is
self-dual: $\chi_\lambda\circ\theta\cong\chi_\lambda$ for every signature $\lambda$.

\subsection{The upper half, one rung at a time}\label{sec:6.3}

For a signature $\lambda$ and a level $m$ on its band we write $u^{(\lambda)}_m$ for the
eigenvalue of $T^{(m)}_{m+1}$ on $W^{(m)}_\lambda$: the rung $m\to m+1$, read at level $m$. As
for $w^{(\lambda)}_n$, which by Lemma~\ref{thm:w-down} is the eigenvalue of $D^{(n)}_s$ on
$W^{(n)}_\lambda$, the subscript is the level acted on. Proposition~\ref{thm:6.9} expresses $w^{(\lambda)}_n$ as a product of these
smaller-step eigenvalues $u^{(\lambda)}_m$, one for each rung $m\to m+1$ between $s$ and
$n$, so that computing $w^{(\lambda)}_n$ comes down to computing one rung at a time.

\begin{proposition}[the ladder of one-step eigenvalues]\label{thm:6.9}
For every signature $\lambda$ and every level $m$ on its band,
\[
\bigl(S^{\,m+1-s}_{1}\bigr)^{2}\,w^{(\lambda)}_{m+1}
\;=\;u^{(\lambda)}_m\;w^{(\lambda)}_{m},
\qquad\text{hence}\qquad
w^{(\lambda)}_n\;=\;\gamma_{s,n}^{-2}\prod_{m=s}^{n-1}u^{(\lambda)}_m,
\]
where $\gamma_{s,n}=\prod_{i=2}^{n-s}S^{\,i}_{1}$ is the number of complete flags
between an $s$-space and an $n$-space containing it. Consequently, for $\lambda\ne\varnothing$, with $s$
parts, the formula of Theorem~\ref{thm:6.15} for $w^{(\lambda)}_n$ is equivalent to
\begin{equation}\label{eq:6.3}
u^{(\lambda)}_m\;=\;[m+1-s]_q[d-s-m]_q\;
q^{\,(r-1)(d-2s-1)+\sigma(\lambda)}\qquad(s\le m\le d-s).
\end{equation}
\end{proposition}

\begin{proof}
By \eqref{eq:2.2}, $M_{sm}M_{m,m+1}=S^{\,m+1-s}_1M_{s,m+1}$. Applying the
transpose to a unit vector $\varphi\in W^{(s)}_\lambda$ gives
\[
\bigl(S^{\,m+1-s}_1\bigr)^2\,w^{(\lambda)}_{m+1}
=\bigl\lVert M_{m,m+1}^{\mathsf T}\,M_{sm}^{\mathsf T}\varphi\bigr\rVert^2
=u^{(\lambda)}_m\,\bigl\lVert M_{sm}^{\mathsf T}\varphi\bigr\rVert^2
=u^{(\lambda)}_m\,w^{(\lambda)}_{m},
\]
since $M_{sm}^{\mathsf T}\varphi$ lies in $W^{(m)}_\lambda$ and
$\lVert M_{sn}^{\mathsf T}\varphi\rVert^2=w^{(\lambda)}_n$. Iterating from
$w^{(\lambda)}_{s}=1$ gives the product formula. The equivalence follows by simple arithmetic.
\end{proof}

At $m=d-s$ the right-hand side of \eqref{eq:6.3} vanishes, as it should: $M_{d-s,d-s+1}^{\mathsf
T}$ maps the $\chi_\lambda$-isotypic part into $W^{(d-s+1)}_\lambda=0$.

\begin{lemma}\label{thm:6.10}
Let $\lambda$ be a signature and $m$ a level on its band. Then
$D^{(m)}_{m-1}$ acts on $W^{(m)}_\lambda$ by $u^{(\lambda)}_{m-1}$ if $m-1$ also lies on the
band, and by $0$ otherwise.
\end{lemma}

\begin{proof}
This is Lemma~\ref{thm:w-down} with $j=m-1$ and $n=m$.
\end{proof}

\begin{corollary}[the ladder is a palindrome]\label{thm:6.11}
For every signature $\lambda$ with $s$ parts and every $s\le m\le d-s-1$,
$u^{(\lambda)}_m=u^{(\lambda)}_{d-m-1}$.
\end{corollary}

\begin{proof}
By Fact~\ref{thm:6.6}, $T^{(m)}_{m+1}(U,U')$, the number of $(m+1)$-spaces containing
$U$ and $U'$, equals the number of $(d-m-1)$-spaces contained in
$U^\perp\cap U'^\perp$, i.e. $D^{(d-m)}_{d-m-1}(U^\perp,U'^\perp)$. Thus
$T^{(m)}_{m+1}\Phi=\Phi D^{(d-m)}_{d-m-1}$ on $\ell^2(\mathcal{X}_{d-m})$. By
Proposition~\ref{thm:6.8}, $\Phi$ maps $W^{(d-m)}_\lambda$ onto
$W^{(m)}_\lambda$, so $u^{(\lambda)}_m$ is the eigenvalue of $D^{(d-m)}_{d-m-1}$ on
$W^{(d-m)}_\lambda$, and this is $u^{(\lambda)}_{d-m-1}$ by Lemma~\ref{thm:6.10}, as the level $d-m-1$
lies on the band of $\lambda$.
\end{proof}

By Proposition~\ref{thm:6.4} and Proposition~\ref{thm:6.9}, the rungs $m\to m+1$ with
$m+1\le d/2$ are known, and Corollary~\ref{thm:6.11} reflects them onto the rungs with
$m\ge d/2$. For even $d$ this covers every rung, the rung $d/2\to d/2+1$ being the mirror
image of $d/2-1\to d/2$. For odd $d$ the rung $m\to m+1$ with $m=(d-1)/2$ is its own
mirror image, and its upper end lies beyond Theorem~\ref{thm:6.2}: Proposition~\ref{thm:6.9}
would need $w^{(\lambda)}_{m+1}$, an eigenvalue on the level $m+1>d/2$.

We compute this rung at its lower end. The operator
$T^{(m)}_{m+1}=M_{m,m+1}M_{m,m+1}^{\mathsf T}$ is an \emph{up-down} operator: its entry at
$(U,U')$ counts the paths $U\subset W\supset U'$ through the level $m+1$, which lies above
$d/2$. The cellular elements $\mathbf c^{(m)}_\mu=M_{\mu m}^{\mathsf T}M_{\mu m}$ are
\emph{down-up} operators: their entry counts the paths $U\supseteq Z\subseteq U'$ through the
submodules of type $\mu$. They never leave the level $m$ and the submodules below it, and
since $2m\le d$, Theorem~\ref{thm:6.2} gives the scalar $A^{(m)}(\mu,\lambda)$ by which each
of them acts on $W^{(m)}_\lambda$. Since $T^{(m)}_{m+1}$ is $G$-equivariant and the cellular
elements form a basis of $\operatorname{End}_G\ell^2(\mathcal{X}_m)$, the up-down operator is a
linear combination of down-up operators,
\[
T^{(m)}_{m+1}=\sum_{\mu\subseteq(r^m)}\beta_\mu\,\mathbf c^{(m)}_\mu,
\qquad\text{hence}\qquad
u^{(\lambda)}_m=\sum_{\mu\subseteq(r^m)}\beta_\mu\,A^{(m)}(\mu,\lambda).
\]
Lemma~\ref{thm:6.12} computes the coefficients $\beta_\mu$. Only $I=\mathbf c^{(m)}_{(r^m)}$ and
the $r$ elements $\mathbf c^{(m)}_{(r^{m-1},j)}$, $0\le j\le r-1$, occur; the latter count the
submodules $Z\cong\mathcal{O}_r^{\,m-1}\oplus\mathcal{O}_j$ of $U\cap U'$, and for $j=0$ this
is $D^{(m)}_{m-1}$ by \eqref{eq:6.1}. Proposition~\ref{thm:6.13} then evaluates the sum.

\begin{lemma}[the one-step operator]\label{thm:6.12}
Let $1\le m\le d/2$ and $z=q^{\,d-2m}$. On $\ell^2(\mathcal{X}_m)$,
\[
T^{(m)}_{m+1}=q^{\,m}[d-2m]_q\,q^{\,(r-1)(d-m-1)}\,I+\mathbf c^{(m)}_{(r^{m-1})}
+(z-1)\sum_{j=1}^{r-1}z^{\,j-1}q^{\,(m-1)j}\,\mathbf c^{(m)}_{(r^{m-1},j)} .
\]
In particular $T^{(m)}_{m+1}=D^{(m)}_{m-1}$ when $d=2m$.
\end{lemma}

\begin{proof}
We compare class values. Let $(U,U')\in \mathcal{X}_m\times \mathcal{X}_m$ have invariant $a$,
in the normal form of Theorem~\ref{thm:3.1}, and $L=\#\{i:a_i\ge1\}$.

\emph{The left-hand side.} An $(m+1)$-space containing $U$ is a line of
$V/U\cong\mathcal{O}_r^{\,d-m}$, and it contains $U'$ if and only if the line contains
$(U+U')/U\cong\bigoplus_{a_i\ge1}\mathfrak{p}^{r-a_i}$ (Proposition~\ref{thm:3.2}). The
submodules of a line form a chain, so $T^{(m)}_{m+1}$ vanishes if $L\ge2$. On the class
$a=(r-k,0^{m-1})$, $0\le k\le r-1$, $(U+U')/U$ is generated by $\pi^{\,k}\bar f_1$ with
$\bar f_1$ primitive, and Lemma~\ref{thm:2.2} gives the value $q^{\,k(d-m-1)}$; on the
diagonal, the class $k=r$, the value is $S^{\,d-m}_1$.

\emph{The cellular elements.} Write $\mathbf c_j=\mathbf c^{(m)}_{(r^{m-1},j)}$ for
$0\le j\le r$, so that $\mathbf c_0=D^{(m)}_{m-1}$ by \eqref{eq:6.1} and
$\mathbf c_r=\mathbf c^{(m)}_{(r^m)}=I$. A module $N$ has $\dim N[\mathfrak{p}]$ parts, where
$N[\mathfrak{p}]=\{x\in N:\pi x=0\}$, and $\dim\pi^{\,r-1}N$ of them equal $r$; the first
does not grow on passing to a submodule. Hence a submodule $Z\subseteq U\cap U'$ with
$|Z|=q^{\,r(m-1)+j}$ is of type $(r^{m-1},j)$ if and only if $\dim\pi^{\,r-1}Z\ge m-1$. As
$U\cap U'\cong\prod_i\mathcal{O}_{r-a_i}$ has $m-L$ parts equal to $r$
(Proposition~\ref{thm:3.2}), every $\mathbf c_j$ vanishes if $L\ge2$. On the class $k$,
$U\cap U'\cong M_k:=F_0\oplus\mathcal{O}_rv$ with $F_0\cong\mathcal{O}_r^{\,m-1}$ and
$\mathcal{O}_rv\cong\mathcal{O}_k$; let $N_j(k)$ be the number of submodules of $M_k$ of type
$(r^{m-1},j)$, the value of $\mathbf c_j$ there. By orders, $N_j(k)=0$ for $j>k$; let
$j\le k$ and $t=k-j$.

If $k<r$, these $Z$ are exactly the kernels of the homomorphisms
$\varphi\colon M_k\to\mathcal{O}_t$ with $\varphi(v)=1$, and $\varphi$ is determined by its
kernel, as the quotient map followed by the isomorphism $M_k/Z\to\mathcal{O}_t$ taking $v$
to $1$. Since $\varphi$ is arbitrary on $F_0$, $N_j(k)=|\mathcal{O}_t|^{m-1}=q^{\,(m-1)t}$.
Indeed, such a kernel has order $q^{\,r(m-1)+j}$ and meets every coset $f+\mathcal{O}_rv$, so
it projects onto $F_0$ along $\mathcal{O}_rv$, and $\pi^{\,r-1}Z$ onto $\pi^{\,r-1}F_0$, of
dimension $m-1$. Conversely, if $Z$ is of type $(r^{m-1},j)$, then
$\pi^{\,r-1}Z=\pi^{\,r-1}M_k=\pi^{\,r-1}F_0$, both having dimension $m-1$; so the projection
of $Z$ to $F_0$ has $m-1$ parts equal to $r$ and is all of $F_0$. Thus
$Z+\mathcal{O}_rv=M_k$, and $M_k/Z$ is cyclic of order $q^{\,t}$, generated by the image of $v$.

If $k=r$, then $M_r\cong U$ is free. A $Z$ of type $(r^{m-1},j)$, $j<r$, contains an
$(m-1)$-space $F$, a summand of $U$ with $U/F\cong\mathcal{O}_r$, so $U/Z$ is cyclic of
order $q^{\,r-j}$, and $Z$ lies in exactly one submodule $W\subset U$ of index $q$. These
$W$ are the preimages of the $[m]_q$ hyperplanes of $U/\pi U\cong\pi^{\,r-1}U$, so
$\dim\pi^{\,r-1}W=m-1$ and $W\cong M_{r-1}$; hence $N_j(r)=[m]_q\,N_j(r-1)$. Altogether
\[
N_j(k)=q^{\,(m-1)(k-j)}\quad(0\le j\le k\le r-1),\qquad
N_j(r)=[m]_q\,q^{\,(m-1)(r-1-j)}\quad(0\le j\le r-1).
\]

\emph{The coefficients, one at a time.} We look for $\beta_0,\dots,\beta_r$ with
$T^{(m)}_{m+1}=\sum_{j=0}^{r}\beta_j\mathbf c_j$, comparing class values. Both sides vanish
on the classes with $L\ge2$, and on the class $k$ only $\mathbf c_0,\dots,\mathbf c_k$ are
nonzero, $\mathbf c_k$ with value $1$; so the class $k$ determines $\beta_k$ once
$\beta_0,\dots,\beta_{k-1}$ are known. Put $x=q^{\,m-1}$ and $y=q^{\,d-m-1}$, so that
$z=y/x$; by the first part, $T^{(m)}_{m+1}$ has value $y^{k}$ on the class $k\le r-1$, and
$S^{\,d-m}_1=[d-m]_q\,y^{\,r-1}$ on the diagonal.

For $k=0$, $U\cap U'$ is an $(m-1)$-space, only $\mathbf c_0$ is nonzero there, with value
$1$, and $T^{(m)}_{m+1}$ has value $1$: so $\beta_0=1$. For $k=1$ the class reads
$\beta_0x+\beta_1=y$, so $\beta_1=y-x$. In general, for $0\le k\le r-1$ the class $k$ reads
\begin{equation}\label{eq:6.12a}
\sum_{j=0}^{k}\beta_j\,x^{\,k-j}=y^{\,k}.
\end{equation}
For $k\ge1$ its terms with $j<k$ are $x$ times the left-hand side of \eqref{eq:6.12a} for
$k-1$, that is $x\,y^{\,k-1}$; hence
\[
\beta_k=y^{\,k}-x\,y^{\,k-1}=(z-1)\,z^{\,k-1}q^{\,(m-1)k}\qquad(1\le k\le r-1).
\]
Finally, the diagonal reads $\beta_r+[m]_q\sum_{j=0}^{r-1}\beta_j\,x^{\,r-1-j}=[d-m]_q\,y^{\,r-1}$,
and the sum is the left-hand side of \eqref{eq:6.12a} for $k=r-1$, that is $y^{\,r-1}$. Hence
\[
\beta_r=\bigl([d-m]_q-[m]_q\bigr)y^{\,r-1}=q^{\,m}[d-2m]_q\,q^{\,(r-1)(d-m-1)} .
\]
\end{proof}

\begin{proposition}[one rung of the lower half]\label{thm:6.13}
Let $\lambda\ne\varnothing$ be a signature with $s$ parts, and let $s\le m\le d/2$. Then
\[
u^{(\lambda)}_m\;=\;[m+1-s]_q[d-s-m]_q\;q^{\,(r-1)(d-1)+2s-|\lambda|} .
\]
\end{proposition}

For $m+1\le d/2$ this also follows from Propositions~\ref{thm:6.4} and \ref{thm:6.9};
the case needed below is $m=(d-1)/2$ for odd $d$.

\begin{proof}
Write $\ell_i=\lambda'_i$, so that $\ell_1=s$ and $\ell_{r+1}=0$, and
$R_j=\sum_{i>j}\ell_i$. By Theorem~\ref{thm:6.2}, the cellular element
$\mathbf c^{(m)}_{(r^{m-1},j)}$ acts on $W^{(m)}_\lambda$ by
$a_j:=A^{(m)}\bigl((r^{m-1},j),\lambda\bigr)$. For $\mu=(r^{m-1},j)$ we have $\mu'_i=m$
for $i\le j$ and $\mu'_i=m-1$ for $i>j$, while $\phi'_i=m$ and $\Phi'_i=d$; so
$\phi'-\mu'$ is the indicator vector of $\{j+1,\dots,r\}$, and the exponent in
Theorem~\ref{thm:6.2} is $\sum_{i>j}(d-1-\ell_i)=(r-j)(d-1)-R_j$. If $j\ge1$, then
$\phi'_1=\mu'_1$, so the two middle factors equal $1$, and $\mu'_i-\mu'_{i+1}$ vanishes
except at $i=j$, where it is $1$, and at $i=r$, where the factor is $\genfrac{[}{]}{0pt}{}{m-1}{m-1}=1$; the
factor at $i=j$ is $\genfrac{[}{]}{0pt}{}{m-\ell_{j+1}}{1}=q^{-(m-1-\ell_{j+1})}[m-\ell_{j+1}]_q$. If $j=0$,
the product equals $1$, the middle factors are $\genfrac{[}{]}{0pt}{}{m-s}{1}\,\genfrac{[}{]}{0pt}{}{d-m-s+1}{1}$, and the
exponent $r(d-1)-|\lambda|$ combines with them. Altogether
\[
\begin{aligned}
a_0&=q^{\,(r-1)(d-1)+2s-|\lambda|}\,[m-s]_q[d-m-s+1]_q,\\
a_j&=q^{\,(r-j)(d-1)-R_j-(m-1-\ell_{j+1})}\,[m-\ell_{j+1}]_q\qquad(1\le j\le r-1).
\end{aligned}
\]
Put $y=q^{\,d-m-1}$, so that $(z-1)z^{\,j-1}q^{\,(m-1)j}=(y-q^{\,m-1})y^{\,j-1}$ and
$q^{\,m}[d-2m]_qq^{\,(r-1)(d-m-1)}=y^{r-1}\bigl([d-m]_q-[m]_q\bigr)$. By
Lemma~\ref{thm:6.12},
\[
u^{(\lambda)}_m
=y^{r-1}\bigl([d-m]_q-[m]_q\bigr)+a_0+\bigl(y-q^{\,m-1}\bigr)\sum_{j=1}^{r-1}
y^{\,j-1}a_j .
\]
Here $[d-m]_q-[m]_q=q^{\,m}[d-2m]_q$ and $y-q^{\,m-1}=(q-1)q^{\,m-1}[d-2m]_q$.
Putting $z_j=q^{\,(1-j)m-R_{j+1}}$ and using $R_j=\ell_{j+1}+R_{j+1}$, the formula for
$a_j$ gives
\[
(q-1)\,q^{\,m-1}y^{\,j-1}a_j=q^{\,(r-1)(d-1)}\bigl(z_{j-1}-z_j\bigr)\qquad(1\le j\le r-1),
\]
and $q^{\,m}y^{r-1}=q^{\,(r-1)(d-1)}z_{r-1}$. The sum telescopes, and since
$z_0=q^{\,m-R_1}=q^{\,m+s-|\lambda|}$,
\[
\begin{aligned}
u^{(\lambda)}_m&=a_0+[d-2m]_q\,q^{\,(r-1)(d-1)}z_0\\
&=q^{\,(r-1)(d-1)+2s-|\lambda|}\Bigl([m-s]_q[d-m-s+1]_q+q^{\,m-s}[d-2m]_q\Bigr).
\end{aligned}
\]
With $a=m-s$ and $b=d-m-s$ the bracket is $[a]_q[b+1]_q+q^{a}[b-a]_q=[a+1]_q[b]_q$, because
$[b+1]_q=[b]_q+q^{b}$, $[a+1]_q=[a]_q+q^{a}$ and $[b]_q-[b-a]_q=q^{\,b-a}[a]_q$.
\end{proof}

Note that $(r-1)(d-1)+2s-|\lambda|=(r-1)(d-2s-1)+\sigma(\lambda)$, so
Proposition~\ref{thm:6.13} is \eqref{eq:6.3} on the lower half of the band. Its ingredient
$a_0$ is the eigenvalue of $D^{(m)}_{m-1}=\mathbf c^{(m)}_{(r^{m-1})}$, i.e. the rung
$m-1\to m$, in accordance with Lemma~\ref{thm:6.10}.

\begin{theorem}[the one-step eigenvalues]\label{thm:6.14}
For every signature $\lambda\ne\varnothing$, with $s$ parts, and every level $s\le m\le d-s$,
\[
u^{(\lambda)}_m\;=\;[m+1-s]_q[d-s-m]_q\;q^{\,(r-1)(d-1)+2s-|\lambda|}.
\]
In particular \eqref{eq:6.3} holds, and $u^{(\lambda)}_m>0$ for $s\le m\le d-s-1$.
\end{theorem}

\begin{proof}
For $m=d-s$ both sides vanish. For $m+1\le d/2$, Proposition~\ref{thm:6.4} gives
$w^{(\lambda)}_m$ and $w^{(\lambda)}_{m+1}$, and Proposition~\ref{thm:6.9} turns them into
\eqref{eq:6.3}. For odd $d$ and $m=(d-1)/2$ this is Proposition~\ref{thm:6.13}. Every
other $m\le d-s-1$ has $m\ge d/2$, and then $m^*=d-m-1$ satisfies $s\le m^*$ and
$m^*+1\le d/2$; by Corollary~\ref{thm:6.11} and the first case,
\[
u^{(\lambda)}_m=u^{(\lambda)}_{m^*}
=[d-m-s]_q[m+1-s]_q\,q^{\,(r-1)(d-1)+2s-|\lambda|} . \qedhere
\]
\end{proof}

\subsection{The blocks}\label{sec:6.4}

\begin{theorem}[the blocks]\label{thm:6.15}
Let $\lambda$ be a signature, $s=s(\lambda)$, and let $n$
be a level on its band. Then
\[
w^{(\varnothing)}_n=S^{\,d-1}_{\,n-1}\,S^{\,n}_{1},
\qquad\qquad
w^{(\lambda)}_n\;=\;S^{\,d-2s}_{\,n-s}\;q^{\,(n-s)\,\sigma(\lambda)}\quad(\lambda\ne\varnothing).
\]
In particular $B^{(\varnothing)}$ is the quotient matrix of the level partition, with
entries $\sqrt{S^{\,d-m}_{\,n-m}S^{\,n}_{\,m}}$; its Perron root is the spectral
radius of $A$, which is a simple eigenvalue, and $\mathbb{P}^{d-1}_{fr}(\mathcal{O}_r)$ is connected.
For a signature $(t)$ with one part, $w^{(t)}_n=S^{\,d-2}_{\,n-1}q^{\,(n-1)(2r-t)}$.
\end{theorem}

\begin{remark}\label{thm:6.16}
The spectral radius of $A$ is the Perron root of $B^{(\varnothing)}$,
whose entries $\sqrt{S^{\,d-m}_{\,n-m}S^{\,n}_{\,m}}$ are the geometric means of the
two degrees between the levels $m<n$ (Section~\ref{sec:2}).
\end{remark}

\begin{remark}[connectivity]\label{thm:6.17}
$\mathbb{P}^{d-1}_{fr}(\mathcal{O}_r)$ is connected for $d\ge3$:
every vertex is adjacent to a line inside it, and two distinct lines $L,L'$ are both
adjacent to the smallest summand containing $L+L'$, which has rank $2\le d-1$. The
level partition is equitable, so the positive Perron eigenvector of
$B^{(\varnothing)}$ lifts to a positive eigenvector of $A$; for a connected graph this
eigenvector belongs to the spectral radius, which is simple by the Perron--Frobenius
theorem. This uses neither Theorem~\ref{thm:4.4} nor Section~\ref{sec:3}.
\end{remark}

\begin{proof}[Proof of Theorem~\ref{thm:6.15}]
For $\lambda=\varnothing$, $\chi_\varnothing$ is the trivial representation and
$W^{(1)}_\varnothing=\mathbb{C}\mathbf 1$; each line lies in $S^{\,d-1}_{\,n-1}$
$n$-spaces, each containing $S^{\,n}_1$ lines, so
$T^{(1)}_n\mathbf 1=S^{\,d-1}_{\,n-1}S^{\,n}_1\mathbf 1$. The assertions on the spectral
radius and connectivity are Remark~\ref{thm:6.17}. Let $\lambda\ne\varnothing$. For $n\le d/2$ the formula is
Proposition~\ref{thm:6.4}. On the whole band, by Theorem~\ref{thm:6.14} the one-step
eigenvalues of $\chi_\lambda$ are given by \eqref{eq:6.3}, and by Proposition~\ref{thm:6.9} this is
equivalent to $w^{(\lambda)}_n=S^{\,d-2s}_{\,n-s}q^{\,(n-s)\sigma(\lambda)}$.
\end{proof}

\begin{proof}[Completion of the proof of Theorem~\ref{thm:4.4}]
By Theorem~\ref{thm:6.15} every $w^{(\lambda)}_n$ is positive on the band,
i.e. condition (P) holds for every signature, and Section~\ref{sec:4} applies. The argument is
not circular: Sections~\ref{sec:6.1}--\ref{sec:6.3} use Section~\ref{sec:4} only through Lemma~\ref{thm:4.1}, which does not depend
on (P).
\end{proof}

\subsection{The universal ladder}\label{sec:6.5}

For $1\le s\le d/2$ and $\sigma\ge0$ put $N=d-2s+1$, and let $L(s,\sigma)=L^{q,r}(s,\sigma)$,
the \textbf{universal ladder}, be the $N\times N$ symmetric matrix with rows and columns
indexed by $0,\dots,N-1$, zero diagonal, and
\begin{equation}\label{eq:L}
L^{q,r}(s,\sigma)_{ij}\;=\;S^{\,j}_{\,i}\,\sqrt{\frac{S^{\,N-1}_{\,j}\,q^{\,j\sigma}}
{S^{\,N-1}_{\,i}\,q^{\,i\sigma}}}\qquad(i<j).
\end{equation}
Since only the
ratios $w_j/w_i$ of $w_j=S^{\,N-1}_{\,j}q^{\,j\sigma}$ enter, $L(d/2,\sigma)=(0)$ for every
$\sigma$, and expanding $S$,
\begin{equation}\label{eq:6.5}
L(s,\sigma)_{ij}\;=\;\binom{j}{i}_q\sqrt{\frac{\binom{N-1}{j}_q}{\binom{N-1}{i}_q}}\;\;
q^{\frac{(j-i)\left[\,\sigma+(r-1)(N-1-j+i)\,\right]}{2}}\qquad(i<j).
\end{equation}

\begin{corollary}[the spectrum in one display]\label{thm:6.18}
For every signature
$\lambda\neq\varnothing$, $B^{(\lambda)}=L\bigl(s(\lambda),\sigma(\lambda)\bigr)$, and
consequently
\[
\operatorname{Spec}(A)\;=\;\operatorname{Spec}\bigl(B^{(\varnothing)}\bigr)
\;\uplus\;\biguplus_{\lambda\neq\varnothing}\;\dim(\chi_\lambda)\cdot
\operatorname{Spec}\,L\bigl(s(\lambda),\sigma(\lambda)\bigr).
\]
In particular a block depends on $\lambda$ only through $s(\lambda)$ and
$\sigma(\lambda)$, and the kernel signatures ($2s=d$) are the case $N=1$, where
$L(d/2,\sigma)=(0)$.
\end{corollary}

\begin{proof}[Proof of Corollary~\ref{thm:6.18}]
Let $\lambda\ne\varnothing$ have $s$ parts,
and index the band by $j=n-s\in\{0,\dots,d-2s\}$. Then the coefficient
$S^{\,n-s}_{\,m-s}$ of Theorem~\ref{thm:4.4} is $S^{\,j}_{\,i}$ with $i=m-s$, $j=n-s$, and by
Theorem~\ref{thm:6.15}, $w^{(\lambda)}_n=S^{\,N-1}_{\,j}q^{\,j\sigma}$ with $N-1=d-2s$. These are
the entries \eqref{eq:L}.
\end{proof}

In particular, two signatures with the same number of parts and the same size have the
same block. The smallest example is $\lambda=(3,1)$ and $\lambda=(2,2)$ at $r=3$, both
with $s=2$ and $\sigma=8$.

\begin{remark}\label{thm:6.19}
The kernel signatures ($2s=d$) have $N=1$, so they appear in the
union of Corollary~\ref{thm:6.18} on the same footing as all other signatures. Comparing with
Theorem~\ref{thm:4.4},
\[
\sum_{s(\lambda)=d/2}\dim\chi_\lambda\;=\;|\mathcal{X}_{d/2}|-|\mathcal{X}_{d/2-1}| .
\]
\end{remark}

\begin{remark}\label{thm:6.20}
The block $B^{(\varnothing)}$, of size $d-1$, is not of the form $L(1,\sigma)$
for any integer $\sigma$. Indeed, with $N=d-1$ and $j=n-1$,
\[
\frac{w^{(\varnothing)}_n}{S^{\,N-1}_{\,j}}
\;=\;\frac{[d-1]_q}{[d-1-j]_q}\;[j+1]_q\;q^{\,2(r-1)j},
\]
and matching the ratio at $j=1$ requires
$q^{\,\sigma}=\frac{[2]_q\,[d-1]_q}{[d-2]_q}q^{\,2(r-1)}$; but the fraction is congruent to
$1$ modulo $q$ and larger than $1$, hence not a power of $q$. For $d\ge5$ not even a
real $\sigma$ satisfies the conditions at $j=1$ and $j=2$ simultaneously.
\end{remark}

For example: for $d=5$ and $s=2$ the block $L(2,\sigma)$ is $2\times2$, with eigenvalues
$\pm q^{\sigma/2}$. For $d=6$ and $s=2$ the block $L(2,\sigma)$ is $3\times3$, and with
$c=q^{\sigma}$ and $b=S^{\,2}_1$ it is
\[
B=\begin{pmatrix}0&\sqrt{bc}&c\\ \sqrt{bc}&0&\sqrt{bc}\\ c&\sqrt{bc}&0\end{pmatrix},
\qquad
\det(xI-B)=(x+c)\,(x^2-cx-2bc),
\]
so its eigenvalues are $-c$ and the roots of $x^2-cx-2bc$. For $r=1$ every signature is $(1^s)$, with
$\sigma=s$.

\section{Subgraphs induced on sets of levels}\label{sec:7}

\begin{theorem}[level-set subgraphs]\label{thm:7.1}
For $\varnothing\ne S\subseteq\{1,\dots,d-1\}$,
the adjacency matrix $A_S$ of $\mathcal{X}_S$ satisfies
\[
\operatorname{Spec}(A_S)\;=\;\biguplus_{\lambda}\;\dim(\chi_\lambda)\cdot
\operatorname{Spec}\Bigl(B^{(\lambda)}\big|_{\,S\cap\,\mathrm{band}(\lambda)}\Bigr),
\]
the union over the signatures whose band meets $S$, where $B|_T$ denotes the
principal submatrix on the rows and columns indexed by $T$.
\end{theorem}

\begin{proof}[Proof of Theorem~\ref{thm:7.1}]
Let $P_S$ be the orthogonal projection of $\ell^2(\mathcal{X})$ onto
$\bigoplus_{n\in S}\ell^2(\mathcal{X}_n)$, so that $A_S$ is $P_SAP_S$ restricted to the range of
$P_S$. Since $P_S$ is $G$-equivariant, it preserves the decomposition of the proof of
Theorem~\ref{thm:4.4} (Section~\ref{sec:4}). In each $A$-invariant subspace $U_i$ there, with orthonormal
basis $(\psi^i_n)_{n\in\mathrm{band}(\lambda)}$ consisting of vectors supported on single
levels, $P_S$ is the coordinate projection onto $\{\psi^i_n:n\in S\}$, and the
compression of $B^{(\lambda)}$ to these coordinates is the principal submatrix
$B^{(\lambda)}|_{S\cap\mathrm{band}(\lambda)}$. The total dimension is
$\sum_\lambda\dim\chi_\lambda\,|S\cap\mathrm{band}(\lambda)|=\sum_{n\in S}|\mathcal{X}_n|$ by
\eqref{eq:3.1}, so this accounts for the whole spectrum of $A_S$.
\end{proof}

\begin{corollary}[two-level subgraphs]\label{thm:7.2}
For $1\le m<n\le d-1$,
\[
\operatorname{Spec}\bigl(\mathcal{X}_{m,n}\bigr)
=\biguplus_{\lambda}\ \dim\chi_\lambda\cdot\bigl\{\pm B^{(\lambda)}_{mn}\bigr\}
\ \uplus\ \{0\}^{\,|\mathcal{X}_m|+|\mathcal{X}_n|-2\sum_\lambda\dim\chi_\lambda},
\]
the unions over the signatures with $s(\lambda)\le m$ and $n\le d-s(\lambda)$, where
\[
B^{(\varnothing)}_{mn}=\sqrt{S^{\,d-m}_{\,n-m}\,S^{\,n}_{\,m}},\qquad
B^{(\lambda)}_{mn}\;=\;S^{\,n-s}_{\,m-s}\,
\sqrt{\frac{S^{\,d-2s}_{\,n-s}}{S^{\,d-2s}_{\,m-s}}}\;\;
q^{\,(n-m)\,\sigma(\lambda)/2}\quad(\lambda\ne\varnothing,\ s=s(\lambda)).
\]
\end{corollary}

In particular, $B^{(\lambda)}_{mn}$ is the singular value of $M_{mn}$ on the
$\chi_\lambda$-isotypic parts, and $\operatorname{rank}M_{mn}=\sum\dim\chi_\lambda$ over
the signatures living on both levels.

\begin{corollary}\label{thm:7.3}
$\chi_\lambda$ contributes to $\operatorname{Spec}(A_S)$ if and only
if $S$ meets the band $[s,d-s]$ of $\lambda$, and it contributes a nonzero eigenvalue
if and only if $S$ meets the band at least twice. In particular, the level $1$ lies
only in the bands of the signatures with $s=1$, so that for $S=\{1,n\}$ the
signatures with $s\ge2$ contribute only zeros; this is why the spectrum of
$\mathcal{X}_{1,n}$ could be computed without them \cite{LP25}.
\end{corollary}

\begin{remark}[{a misprint in \cite[Thm. 5.1]{LP25}}]\label{thm:7.4}
For $s=1$, the spectrum of the
bipartite graph $\mathcal{X}_{1,n}$ is $\pm\sqrt{w^{(\lambda)}_n}$ over the signatures
with $s=1$, together with the kernel (Corollary~\ref{thm:7.2}), so Theorem~\ref{thm:6.15} for $s=1$ is
equivalent to \cite[Thm. 5.1]{LP25}. Every nonzero eigenvalue printed there (in the notation
of \cite{LP25}, which is ours) is the true one multiplied by $q^{E/2}$, where
\[
E=(r-1)(d-n)(d-n-1).
\]
This is visible already in the first row: a connected biregular bipartite graph has
largest eigenvalue $\sqrt{\deg L\cdot\deg U}$, and with the degrees
$S^{\,d-1}_{\,n-1}$ and $S^{\,n}_1$ recorded in \cite{LP25} this has $q$-exponent
$(r-1)(n-1)(d-n+1)$, whereas the printed exponent
$(r-1)\bigl((d-1)^2-(d-1)n+2n-2\bigr)$ exceeds it by $E$. The correct eigenvalues are
\[
\pm\sqrt{S^{\,d-1}_{\,n-1}S^{\,n}_{1}}
\qquad\text{and}\qquad
\pm\sqrt{S^{\,d-2}_{\,n-1}\;q^{(n-1)(2r-t)}},\quad t=1,\dots,r,
\]
and the multiplicities printed there are correct. Since $E=0$ when $r=1$ or $n=d-1$,
and $E$ is the same for all nonzero rows, the misprint affects neither [LP25,
Cor. 5.2], nor the spectral expansion computed there, nor the isospectrality results.
\end{remark}

\begin{remark}[dependence on the ring]\label{thm:7.5}
Of the ring $\mathcal{O}_r$, Sections~\ref{sec:3}--\ref{sec:4}
use only that every nonzero element is a unit times $\pi^{\,i}$ for a unique
$i<r$, with $\pi^{\,r}=0$ (Fact~\ref{thm:2.1}(4)), and that the residue field has $q$
elements; Section~\ref{sec:dims} and Sections~\ref{sec:6.1}--\ref{sec:6.3} use in addition,
through \cite{OS06} and Theorem~\ref{thm:6.2}, that $\mathcal{O}_r\cong\mathfrak{O}/\mathfrak{m}^r$ for the ring of integers
$\mathfrak{O}$ of a non-archimedean local field (Section~\ref{sec:2}), which is the generality
in which \cite{BO12} and \cite{OS06} are proved; and Fact~\ref{thm:6.6} holds over any commutative ring.
The answers of Theorems~\ref{thm:4.4}--\ref{thm:7.1} and Corollaries~\ref{thm:6.18} and \ref{thm:7.2} are expressed in $q$,
$r$, $d$ and $\lambda$ only, so the spectrum of every $\mathcal{X}_S(\mathcal{O}_r)$
depends on the ring only through $q$ and $r$. For instance $\mathbb{Z}/2^r$,
$\mathbb{F}_2[t]/(t^r)$ and $\mathbb{Z}[i]/(1+i)^r$ all give the same spectra.
\end{remark}

\begin{theorem}[isospectral non-isomorphic geometries]\label{thm:7.6}
Let $d\ge3$.
\begin{enumerate}
\item If $\mathcal{O}_r=\mathcal{O}/\mathfrak{p}^r$ and
  $\mathcal{O}'_r=\mathcal{O}'/\mathfrak{p}'^{\,r}$ are as above, with the same
  residue order $q$ and the same $r$, then for every nonempty set $S$ of levels the
  graphs $\mathcal{X}_S(\mathcal{O}_r)$ and $\mathcal{X}_S(\mathcal{O}'_r)$ are
  isospectral.
\item Let $p$ be a prime and $r\ge2$. For every set $S$ of at least two levels the graphs
  \[
  \mathcal{X}_S\bigl(\mathbb{Z}/p^r\bigr)\qquad\text{and}\qquad
  \mathcal{X}_S\bigl(\mathbb{F}_p[t]/(t^r)\bigr)
  \]
  are isospectral and non-isomorphic. In particular this holds for the two-level graphs
  $\mathbb{P}^{\,d-1}_{m,n}$, $1\le m<n\le d-1$, and for the full geometries
  $\mathbb{P}^{\,d-1}_{fr}(\mathbb{Z}/p^r)$ and $\mathbb{P}^{\,d-1}_{fr}(\mathbb{F}_p[t]/(t^r))$.
\end{enumerate}
\end{theorem}

\begin{proof}[Proof of Theorem~\ref{thm:7.6}]
Part (1) follows from Theorem~\ref{thm:7.1} and Remark~\ref{thm:7.5}. In part
(2), isospectrality is the case $\mathcal{O}=\mathbb{Z}$, $\mathfrak{p}=(p)$,
$\mathcal{O}'=\mathbb{F}_p[t]$, $\mathfrak{p}'=(t)$ of part (1), with $q=p$.
Non-isomorphism is the argument of \cite[\S6]{LP25} (see the discussion following Thm. 6.1
there).
\end{proof}


\begin{thebibliography}{BO07}

\bibitem[BO07]{BO07} U. Bader and U. Onn, \emph{Geometric representations of $\mathrm{GL}(n,R)$, cellular
Hecke algebras and the embedding problem}, J. Pure Appl. Algebra \textbf{208} (2007),
905--922.

\bibitem[BO12]{BO12} U. Bader and U. Onn, \emph{On some geometric representations of
$\mathrm{GL}_N(\mathfrak{O})$}, Comm. Algebra \textbf{40} (2012), 3169--3191;
arXiv:math/0404408.

\bibitem[Gar73]{Gar73} H. Garland, \emph{$p$-adic curvature and the cohomology of discrete subgroups of
$p$-adic groups}, Ann. of Math. \textbf{97} (1973), 375--423.

\bibitem[Kan72]{Kan72} W. M. Kantor, \emph{On incidence matrices of finite projective and affine spaces},
Math. Z. \textbf{124} (1972), 315--318.

\bibitem[KP22]{KP22} T. Kaufman and O. Parzanchevski, \emph{Free flags over local rings and powering of
high-dimensional expanders}, Int. Math. Res. Not. IMRN \textbf{2022}, no. 19, 14741--14769;
arXiv:1909.02473.

\bibitem[Lew24]{Lew24} A. Lew, \emph{Asymptotic behavior of Laplacian eigenvalues of subspace inclusion
graphs}, J. London Math. Soc. \textbf{109} (2024); arXiv:2308.08397.

\bibitem[LP25]{LP25} Y. Lavi and O. Parzanchevski, \emph{Geometry over finite local rings: rigidity and
isospectrality}, arXiv:2510.09324.

\bibitem[Opp18]{Opp18} I. Oppenheim, \emph{Local spectral expansion approach to high dimensional
expanders, Part I: Descent of spectral gaps}, Discrete Comput. Geom. \textbf{59} (2018),
293--330.

\bibitem[OS06]{OS06} U. Onn and J. V. Stokman, \emph{Quantum dimensions and their non-Archimedean
degenerations}, Int. Math. Res. Pap. IMRP \textbf{2006}, Art. ID 54701.

\bibitem[Vol25]{Vol25} P. Vollrath, \emph{Full resolution to Papikian's conjecture}, arXiv:2510.15161.

\end{thebibliography}
\end{document}